 %&amstex                    
\input amstex\documentstyle{amsppt}  
\pagewidth{12.5cm}\pageheight{19cm}\magnification\magstep1
\topmatter
\title Elliptic elements in a Weyl group: a homogeneity property\endtitle
\author G. Lusztig\endauthor
\address{Department of Mathematics, M.I.T., Cambridge, MA 02139}\endaddress
\thanks{Supported in part by the National Science Foundation}\endthanks
\endtopmatter   
\document

\define\ul{\un l}

\define\uuG{\un{\un G}}

\define\uWW{\un\WW}

\define\dsv{\dashv}

\define\pe{\perp}

\define\qua{\quad}

\define\lb{\linebreak}

\define\eSb{\endSb}

\define\bin{\binom}
\define\op{\oplus}
\define\h{\frac{\hphantom{aa}}{\hphantom{aa}}}   

\define\part{\partial}
\define\em{\emptyset}

\define\m{\mapsto}
\define\do{\dots}

\define\sub{\subset}    

\define\T{\times}
\define\ti{\tilde}
\define\nl{\newline}
\redefine\i{^{-1}}

\define\un{\underline}

\define\sg{\text{\rm sgn}}

\define\a{\alpha}
\redefine\b{\beta}

\define\g{\gamma}
\redefine\d{\delta}
\define\e{\epsilon}

\define\p{\pi}

\define\ps{\psi}

\define\s{\sigma}

\define\k{\kappa}
\redefine\l{\lambda}
\define\z{\zeta}
\define\x{\xi}

\define\Ph{\Phi}

\define\kk{\bold k}

\define\nn{\bold n}

\define\CC{\bold C}

\define\FF{\bold F}

\define\NN{\bold N}

\define\WW{\bold W}
\define\ZZ{\bold Z}

\define\cb{\Cal B}
\define\cc{\Cal C}

\define\cf{\Cal F}

\define\ci{\Cal I}

\define\cm{\Cal M}

\define\co{\Cal O}
\define\cp{\Cal P}

\define\car{\Cal R}

\define\cw{\Cal W}
\define\cz{\Cal Z}

\define\fB{\frak B}

\define\ta{\ti a}

\define\tg{\ti g}

\define\tv{\ti v}
\define\tw{\ti w}

\define\tA{\ti A}
\define\tB{\ti B}

\define\tN{\ti N}

\define\GE{Ge}
\define\GH{GH}
\define\GP{GP}
\define\HS{HS}
\define\IC{L1}
\define\CS{L2}
\define\GF{L3}
\define\WE{L4}
\define\LUE{L\"u}
\define\SH{Sh}

\head Introduction\endhead
\subhead 0.1\endsubhead
Let $G$ be a connected reductive algebraic group over an algebraically closed field $\kk$ of characteristic $p$.
Let $\WW$ be the Weyl group of $G$. Let $\cb$ the variety of Borel subgroups of $G$. For each $w\in\WW$ let 
$\co_w$ be the corresponding $G$-orbit in $\cb\T\cb$. Let $\ul:\WW@>>>\NN$ be the standard length function. Let 
$\uWW$ be the set of conjugacy classes in $\WW$. For $C\in\uWW$ let $d_C=\min_{w\in C}\ul(w)$ and let 
$C_{min}=\{w\in C;l(w)=d_C\}$; let $\Ph(C)$ be the unipotent class in $G$ associated to $C$ in \cite{\WE, 4.1}. 
For any conjugacy class $\g$ in $G$ and any $w\in\WW$ we set $\fB^\g_w=\{(g,B)\in\g\T\cb;(B,gBg\i)\in\co_w\}$; 
note that $G$ acts on $\fB^\g_w$ by $x:(g,B)\m(xgx\i,xBx\i)$. 

For $w\in\WW$ let $\mu(w)$ be the dimension of the fixed point space of $w:V@>>>V$ where $V$ is the reflection
represention of the Coxeter group $\WW$. We say that $w$ or its conjugacy class is elliptic if $\mu(w)=0$.
Let $\uWW_{el}$ be the set of elliptic conjugacy classes in $\WW$.

The following is the main result of this paper. 

\proclaim{Theorem 0.2} Let $C\in\uWW_{el}$ and let $w\in C_{min}$, $\g=\Ph(C)$. Then $\fB^\g_w$ is a single 
$G$-orbit. 
\endproclaim
In the case where $p$ is not a bad prime for $G$, the weaker result that $\fB^\g_w$ is a union of finitely 
many $G$-orbits is already known from \cite{\WE, 5.8(a),(b)}. 

\subhead 0.3\endsubhead
In the setup of 0.2 let $g\in\g$ and let $\cb_g^w=\{B\in\cb;(B,gBg\i)\in\co_w\}$. Let $Z(g)$ be the centralizer
of $g$ in $G$. The folowing result is an immediate consequence of 0.2.

(a) {\it $\cb_g^w$ is a single orbit for the conjugation action of $Z(g)$.}
\nl
It is likely that Theorem 0.2 (and its consequence (a)) continues to hold if $C$ is a not necessarily elliptic 
conjugacy class. See 4.2 for a partial result in this direction.

\subhead 0.4\endsubhead
The proof of the theorem is given in 3.2. It is a case by case argument. The proof for classical groups is easy
for type $A$, relatively easy for type $C$ and much more complicated for types $B,D$. The proof for exceptional 
groups can be reduced (using representation theory, as in \cite{\WE}) to a computer calculation. This 
calculation uses the character tables of Hecke algebras available through the CHEVIE package (see \cite{\GH}),
the tables of Green polynomials in good characteristic (see \cite{\LUE}) and the analogous tables (provided to
me by F. L\"ubeck) in bad characteristic; the fact that these last tables, computed using the algorithm in 
\cite{\CS, Ch.24}, give indeed the Green functions, is proved by M. Geck \cite{\GE} using earlier results in 
\cite{\CS,\GF,\SH}. I thank Gongqin Li for her help with programming in GAP3.

\subhead 0.5. Notation\endsubhead 
Let $S=\{w\in\WW;\ul(w)=1\}$. For any subset $K$ of $S$ let $\WW_K$ be the subgroup
of $\WW$ generated by $K$ and let $\cp_K$ be the conjugacy class of parabolic subgroups of $G$ determined by $K$.
(For example, $\cp_\em=\cb$.) For $B\in\cb$ let $P^K_B$ be the unique subgroup in $\cp_K$ that contains $B$.
Let $G_{ad}$ be the adjoint group of $G$. Let $\cz_G$ be the centre of $G$.
For any nilpotent endomorphism $N$ of a finite dimensional vector space $V$ 
we denote by $\cm(N,V)$ the multiset consisting of the sizes of Jordan blocks of $N$. For $n\in\NN$ define 
$\k_n\in\{0,1\}$ by $n-\k_n\in2\NN$. For $i\in\ZZ-\{0\}$ we define $\sg(i)\in\{1,-1\}$ by the condition that 
$\sg(i)i>0$. If $X$ is a finite set and $f:X@>>>X$ is a map, we set $X^f=\{x\in X;f(x)=x\}$.

\head 1. Isometry groups\endhead
\subhead 1.1\endsubhead
Let $V$ be a $\kk$-vector space of finite dimension $\nn\ge3$. We set $\k=\k_\nn$. Let $n=(\nn-\k)/2$. Assume that
$V$ has a fixed bilinear form $(,):V\T V@>>>\kk$ and a fixed quadratic form $Q:V@>>>\kk$ such that (i) or (ii) 
below holds:

(i) $Q=0$, $(x,x)=0$ for all $x\in V$, $V^\pe=0$;

(ii) $Q\ne0$, $(x,y)=Q(x+y)-Q(x)-Q(y)$ for $x,y\in V$, $Q:V^\pe@>>>\kk$ is injective.
\nl
Here, for any subspace $V'$ of $V$ we set $V'{}^\pe=\{x\in V;(x,V')=0\}$. In case (ii) it follows that $V^\pe=0$ 
unless $\k=1$ and $p=2$ in which case $\dim V^\pe=1$. 

An element $g\in GL(V)$ is said to be an isometry if $(gx,gy)=(x,y)$ for all $x,y\in V$ and $Q(gx)=Q(x)$ for all 
$x\in V$. Let $Is(V)$ be the group of all isometries of $V$ (a closed subgroup of $GL(V)$). 

\subhead 1.2\endsubhead
Let $p_1\ge p_2\ge\do\ge p_\s$ (or $p_*$) be a descending sequence of integers $\ge1$ such that
$p_1+p_2+\do+p_\s=n$. If $\k=1$ we set $p_{\s+1}=1/2$.

Let $g\in Is(V)$. A collection of vectors $w^t_i$ ($t\in[1,\s+\k],i\in\ZZ$) in $V$ is said to be 
$(g,p_*)$-adapted if 

(a) $w^t_{i+1}=gw^t_i$ for all $t,i$;

(b) $(w^t_i,w^t_j)=0$ if $|i-j|<p_t$, $(w^t_i,w^t_j)=1$ if $j-i=p_t$ ($t\in[1,\s]$, $i,j\in\ZZ$);

(c) $(w^t_i,w^r_j)=0$ if $0\le i-j+p_r<2p_t$ and $1\le t<r\le\s$;

(d) $(w^{\s+1}_i,w^{\s+1}_i)=2$ for all $i$ (if $\k=1$);

(e) $(w^t_i,w^{\s+1}_j)=0$ if $\k=1$, $0\le i-j<2p_t$ and $1\le t\le\s$.

(f) $Q(w^t_i)=0$ if $t\in[1,\s],i\in\ZZ$ and $Q(w^{\s+1}_i)=1$ if $\k=1,i\in\ZZ$.

\subhead 1.3\endsubhead
We preserve the setup of 1.2. We show:

(a) {\it $\{w^x_i;x\in[1,\s+\k],i\in[0,2p_x-1]\}$ is a basis of $V$.}

(b) {\it Assume that $e,f\in[1,\s+\k]$, ($e\le f$) and that the subspace $\cw_{e,f}$ of $V$ spanned by
$\{w^x_i;x\in[e,f],i\in[0,2p_x-1]\}$ is $g$-stable. Then the radical $\car$ of $(,)|_{\cw_{e,f}}$ is $0$ unless 
$\k=1,f=\s+1,p=2$, in which case $\car=V^\pe$.}
\nl
We prove (b). The proof is similar to that of \cite{\WE, 3.3(iv),(vi)}. Let $[e,f]^*=[e,f]\cap[1,\s]$.
Assume that $c^x_i\in\kk$ $(x\in[e,f]^*,i\in[0,2p_x-1])$ are not all zero and that
$$\sum_{x\in[e,f]^*,i\in[0,2p_x-1]}c^x_i(w^x_i,w^y_j)=0\tag c$$
for any $y\in[e,f]^*,j\in[0,2p_y-1]$. Let 
$$i_0=\min\{i\in\NN;c^r_i\ne0\text{ for some }r\in[e,f]^*\text{ such that }i\le2p_r-1\},$$ 
$$X=\{r\in[e,f]^*;c^r_{i_0}\ne0,i_0\le2p_r-1\}.$$ 
We have $X\ne\em$. Let $r_0$ be the largest number in $X$.
Since $\cw_{e,f}$ is $g$-stable, $w^{r_0}_{i_0+p_{r_0}}$ is a linear combination of elements $w^y_j$,
$y\in[e,f]^*,j\in[0,2p_j-1]$ (which, by (c), have inner product zero with
$\sum_{x\in[e,f]^*,i\in[0,2p_x-1]}c^x_iw^x_i$) and (if $\k=1,f=\s+1$) of $w^{\s+1}_0$ 
(which also has inner product zero with $\sum_{x\in[e,f]^*,i\in[0,2p_x-1]}c^x_iw^x_i$). Hence
$$\sum_{x\in[e,f]^*,i\in[0,2p_x-1]}c^x_i(w^x_i,w^{r_0}_{i_0+p_{r_0}})=0.$$
This can be written as follows:
$$\sum_{r\in X}c^r_{i_0}(w^r_{i_0},w^{r_0}_{i_0+p_{r_0}})+
\sum_{r\in[e,f]^*;i\in[i_0+1,2p_r-1]}c^r_i(w^r_i,w^{r_0}_{i_0+p_{r_0}}).$$
If $r\in X$, $r\ne r_0$, we have $(w^r_{i_0},w^{r_0}_{i_0+p_{r_0}})=0$ (using $r<r_0$). If $r\in[e,f]^*$ 
and $i\in[i_0+1,2p_r-1]$, we have $(w^r_i,w^{r_0}_{i_0+p_{r_0}})=0$: if $r<r_0$, we have $1\le i-i_0\le2p_r-1$; if
$r\in[e,f]^*$, $r\ge r_0$ we have 
$$1\le-2p_r+1+p_{r_0}+p_r\le i_0-i+p_{r_0}+p_r\le-1+p_{r_0}+p_r\le2p_{r_0}-1.$$
We see that
$$0=c^{r_0}_{i_0}(w^{r_0}_{i_0},w^{r_0}_{i_0+p_{r_0}})=c^{r_0}_{i_0};$$
this contradicts $c^{r_0}_{i_0}\ne0$.

If $f\le\s$ the previous argument shows that the symmetric matrix

$(w^x_i,w^y_j)_{x,y\in[e,f],i\in[0,2p_x-1],j\in[0,2p_y-1]}$
\nl
is nonsingular hence (b) holds.

Assume now that $f=\s+1$ so that $\k=1$. Let $\x\in\car$. Then
there exist $c^x_i\in\kk$ $(x\in[e,\s],i\in[0,2p_x-1])$ and $c\in\kk$ are such that
$$\x=\sum_{x\in[e,\s],i\in[0,2p_x-1]}c^x_iw^x_i+cw^{\s+1}_0.$$
We have
$$\sum_{x\in[e,\s],i\in[0,2p_x-1]}c^x_i(w^x_i,w^y_j)+c(w^{\s+1}_0,w^y_j)=0\tag d$$
for any $y\in[e,\s+1],j\in[0,2p_y-1]$.
For $y\in[e,\s],j\in[0,2p_y-1]$, (d) becomes 
$$\sum_{x\in[e,\s],i\in[0,2p_x-1]}c^x_i(w^x_i,w^y_j)=0;$$
this implies by the first part of the argument that $c^x_i=0$ for all $x\in[e,\s],i\in[0,2p_x-1]$.
Thus $\x=cw^{\s+1}_0$ and (d) implies $0=c(w^{\s+1}_0,w^{\s+1}_0)=2c$. If $p\ne2$ this implies $c=0$ so that 
$\x=0$. Thus in this case (b) holds.
If $p=2$ we see that $\car\sub\kk w^{\s+1}_0$; conversely it is clear
that $\kk w^{\s+1}_0\sub\car$ hence (b) holds again. (Note that $w^{\s+1}_0\ne0$ since $Q(w^{\s+1}_0)=1$.)

We prove (a). We use (b) and its proof with $e=1,f=\s$. If $\k=0$, the nonsingularity of the symmetric matrix
$(w^x_i,w^y_j)_{x,y\in[1,\s],i\in[0,2p_x-1],j\in[0,2p_y-1]}$ shows that the vectors in (a) are linearly 
independent hence they form a basis of $V$ (the number of these vectors is $\sum_{r\in[1,\s]}2p_r=\dim V$).
Now assume that $\k=1$. Assume that 
$$\sum_{x\in[e,\s],i\in[0,2p_x-1]}c^x_iw^x_i+cw^{\s+1}_0=0$$
where $c^x_i\in\kk$ $(x\in[e,\s],i\in[0,2p_x-1])$ and $c\in\kk$ are not all zero.
By the proof of (b) we must have $c^x_i=0$ for all $x\in[e,\s],i\in[0,2p_x-1]$. Hence 
$cw^{\s+1}_0=0$. Since $w^{\s+1}_0\ne0$ it follows that $c=0$, a contradiction. We see that the vectors in (a)
are linearly independent; since their number equals $\dim V$, they form a basis. This completes the proof of (a).

\subhead 1.4\endsubhead
We now assume that $g$ is unipotent. We set $N=g-1:V@>>>V$.
Assume that $k>0$ that $d\in[1,\s]$ is such that $2p_d\ge k\ge2p_{d+1}$ (convention: $p_{\s+1}=0$ if $\k=0$) and
that $\dim N^kV=\sum_{r\in[1,d]}(2p_r-k)$. Let $\cw=\cw_{1,d}$, $\cw'=\cw_{d+1,\s+\k}$ (convention: if $\k=0,d=\s$
then $\cw'=0$). We have the following result. 

(a) {\it $\cw,\cw'$ are $g$-stable, $\cw'=\cw^\pe$, $g:\cw@>>>\cw$ has exactly $d$ Jordan blocks.}
\nl
The proof is exactly the same as that of \cite{\WE, 3.5(b)} if we replace $v'_r$ by $w^r_0$.

\subhead 1.5\endsubhead
In the setup of 1.4 we assume that either $Q=0$ or $p=2$. We also assume that $\cm(N,V)$ consists of 
$2p_1\ge2p_2\ge\do\ge2p_\s$ (and $1$ if $\k=1$). Then for any $k\ge0$ we have 
$\dim N^k(V)=\sum_{r\in[1,\s+\k]}\max(2p_r-k,0)$.
For any $r\in[1,\s+\k]$ let $X_r$ be the subspace of $V$ spanned by $\{w^r_i;i\in[0,2p_r-1]\}$. Note that
$V=\op_{r\in[1,\s+\k}X_r$ (see 1.3(a)). We have the following result:

(a) {\it For any $r\in[1,\s+\k]$, $X_r$ is a $g$-stable subspace of $V$ and for any $r\ne t$ in $[1,\s+\k]$ we
have $(X_r,X_t)=0$.}
\nl
This is deduced from 1.4(a) in the same way as \cite{\WE, 3.5(c)} is deduced from \cite{\WE, 3.5(b)}.

We have the following result.

\proclaim{Proposition 1.6} In the setup of 1.5, assume that $w^t_i$ ($t\in[1,\s+\k],i\in\ZZ$) is 
$(g,p_*)$-adapted. Then:

(a) $(w^t_i,w^t_j)=\sg(j-i)\bin{|j-i|+\p-1}{|j-i|-\p}$ if $|j-i|\ge\p$, $(w^t_i,w^t_j)=0$ if $|j-i|<\p$ 
$(t\in[1,\s],\p=p_t)$;    

(b) $(w^{\s+1}_i,w^{\s+1}_j)=0$ if $\k=1$, $i,j\in\ZZ$;

(c) $(w^t_i,w^r_j)=0$ if $t\ne r$ in $[1,\s+\k]$, $i,j\in\ZZ$.
\endproclaim
We prove (a). Let $t,\p$ be as in (a). Now $X_t$ is $g$-stable (see 1.5(a)) of dimension $2\p$ and $N$ acts on
it as a single Jordan block of size $2\p$. Hence $N^{2\p}w^t_0=0$ that is,
$$\sum_{k\in[0,2\p]}n_kw^t_{k}=0\tag d$$
where $n_k=(-1)^k\bin{2\p}{k}$. Aplying $(,w^t_{\p+s})$ (with $s\in\ZZ_{>0}$) to (d) we obtain
$$(w^t_0,w^t_{s+\p})+\sum_{k\in[1,2\p];k\le s}n_k(w^t_0,w^t_{s+\p-k})=0.$$ 
(Note that if $k>s$ then $(w^t_0,w^t_{s+\p-k})=0$.) 
This can be viewed as an inductive formula for $(w^t_0,w^t_{s+\p})$ (for $k\in[1,2\p],k\le s$ we have
$s+\p-k\in[\p,s+\p-1]$). We show that $\bin{2\p+s-1}{s}$ satisfies the same inductive formula hence
$(w^t_0,w^t_{s+\p})=\bin{2\p+s-1}{s}$. It is enough to show that
$$\sum_{k\in[0,2\p];k\le s}(-1)^k\bin{2\p}{k}\bin{2\p+s-k-1}{s-k}=0$$
for $s\ge1$ or, setting $m=s-k$, that
$$\sum_{s\ge0}\sum_{k\in[0,2\p];m\ge0;k+m=s}(-1)^k\bin{2\p}{k}\bin{2\p+m-1}{m}T^s=1$$
where $T$ is an indeterminate. An equivalent statement is
$$(\sum_{k\in[0,2\p]}(-1)^k\bin{2\p}{k}T^k)(\sum_{m\ge0}\bin{2\p+m-1}{m}T^m)=1.$$
This folows from the identity $\sum_{m\ge0}\bin{M+m-1}{m}T^m=(1-T)^{-M}$ (for $M\ge1$) which is easily verified.
This proves (a).

We prove (b). We now have $\k=1$. From 1.5(a) we see that $X_{\s+1}$ is $g$-stable. It is $1$-dimensional 
hence $g$ acts on it as the identity map. Thus $w^{\s+1}_i=w^{\s+1}_0$ and (b) is reduced to the case where
$i=j=0$ where it follows from the definition.

Now (c) follows from 1.5(a). The proposition is proved.

\subhead 2. Orthogonal groups in odd characteristic\endsubhead
\subhead 2.1\endsubhead
In this section we assume that $\kk$ has characteristic $\ne2$.
We fix a map $\kk\m\kk$, $\l\m\sqrt{\l}$ such that $(\sqrt{\l})^2=\l$ for all $\l\in\kk$. We fix $\k\in\{0,1\}$. 
Assume that $p_1\ge p_2\ge\do\ge p_\s$ is a sequence in $\ZZ_{>0}$. When $\k=1$ we set $p_{\s+1}=1/2$. 

As in \cite{\WE, 1.6}, we define a function $\ps:[1,\s]@>>>\{-1,0,1\}$ as follows.

(i) If $t\in[1,\s]$ is odd and $p_t<p_x$ for any $x\in[1,t-1]$ then $\ps(t)=1$;

(ii) if $t\in[1,\s]$ is even and $p_x<p_t$ for any $x\in[t+1,\s]$, then $\ps(t)=-1$;

(iii) for all other $t\in[1,\s]$ we have $\ps(t)=0$.
\nl
For any $y,x\in[1,\s+\k]$ and $i,j\in\ZZ$ we will define $|^y_i:^x_j|\in\kk$ in 2.2-2.8. We require that 
$|^y_i:^x_j|=|^x_j:^y_i|$. Hence it is enough to define $|^y_i:^x_j|$ under the following assumptions (which will
be in force until the end of 2.8):

$\p$ is a fixed number equal to one of the $p_1,p_2,\do,p_{\s+\k}$;

$p_x=\p$, $y\le x$ (and $i\le j$ if $y=x$);

$|^t_i:^r_j|$ is already defined whenever $p_t>\p,p_r>\p$;
\nl
(the last assumption is empty if $\p=p_1$).

In the case where $p_1>\p$ we define integers $a,b\in[1,\s]$ by the following requirements:

$p_b>p_{b+1}=\p$;

$p_a>\p$, $\ps(a)=1$ (hence $a$ is odd) and $a$ is maximal with these properties;
\nl
we have $a\le b$ and we set $I_\p=[a,b]$.

For any $k\in[0,2\p]$ we set $n_k=(-1)^k\bin{2\p}{k}$.

\subhead 2.2 \endsubhead
Assume that $p_y>\p$ and that for some even $r\in[y,x-1]$ we have $p_r>p_{r+1}$. We set $|^y_i:^x_j|=0$ for all
$i,j\in\ZZ$.

\subhead 2.3\endsubhead
Assume that $p_y>\p\ge1$ and that for any even $r\in[y,x-1]$ we have $p_r=p_{r+1}$. Note that $p_1>\p$ hence 
$I_\p=[a,b]$ is defined as in 2.1. We set $I=I_\p$. We have $y\in I$. 

We define some auxiliary elements $\a^r_h\in\kk$ $(r\in I;h\in[0,p_r-\p-1])$, $\b^r_{2p_r-2\p-h}\in\kk$ 
$(r\in I;h\in[1,p_r-\p])$ by induction on $h$, by the following equations:
$$\align&\a^r_h+\sum'_{r'\in I;r'>r;p_r=p_r'}\a^{r'}_h|_0^{r'}:^r_{p_r}|
=-\sum_{k\in[0,2\p]}n_k|_{k+2p_a-2\p}^a:^r_{p_r+h}|\\&
-\sum'_{r'\in I;i\in[0,p_{r'}-\p-1];i<h;k\in[0,2\p]}\a^{r'}_in_k|_{k+i}^{r'}:^r_{p_r+h}|\\&
-\sum''_{r'\in I;j\in[1,p_{r'}-\p];j<h;k\in[0,2\p]}\b^{r'}_{2p_{r'}-2\p-j}n_k|_{k+2p_{r'}-2\p-j}^{r'}:^r_{p_r+h}|
\tag i\endalign$$
for $r\in I,h\in[0,p_r-\p-1]$ and
$$\align&\b^r_{2p_r-2\p-h}+\sum''_{r'\in I;r'<r}\b^{r'}_{2p_{r'}-2\p-h}|_{2p_{r'}}^{r'}:^r_{p_r}|
=-\sum_{k\in[0,2\p]}n_k|_{k+2p_a-2\p}^a:^r_{p_r-h}|\\&
-\sum'_{r'\in I;i\in[0,p_{r'}-\p-1];i<h-1;k\in[0,2\p]}\a^{r'}_in_k|_{k+i}^{r'}:^r_{p_r-h}|\\&
-\sum''_{r'\in I;j\in[1,p_{r'}-\p];j<h;k\in[0,2\p]}\b^{r'}_{2p_{r'}-2\p-j}n_k|_{k+2p_{r'}-2\p-j}^{r'}:^r_{p_r-h}|
\tag ii\endalign$$  
for $r\in I,h\in[1,p_r-\p]$. Note that the right hand sides of (i),(ii) can be assumed to be known from the 
induction hypothesis. Also if $h=0$ the right hand side of (i) is $0$ and for $h=1$ the right hand side of (ii) 
is $0$. Thus we may assume that
$$\a^r_h+\sum_{r'\in I;r'>r;p_r=p_{r'}}|_{p_{r'}}^{r'}:^r_{2p_r}|\a^{r'}_h$$
is known (this determines $\a^r_h$ by descending induction on $r\in I$) and that 
$$\b^r_{2p_r-2\p-h}+\sum_{r'\in I;r'<r}|_{2p_{r'}}^{r'}:^r_{p_r}|\b^{r'}_{2p_{r'}-2\p-h}$$
is known (this determines $\b^r_{2p_r-2\p-h}$ by induction on $r\in I$).

We now define $\ta^t_j\in\kk$ for $t\in I,j\in[0,2p_t-2\p-1]$ as follows. If $j\in[0,p_t-\p-1]$ then 
$\ta^t_j=\a^t_j$. If $j\in[p_t-\p,2p_t-2\p-1]$ then $\ta^t_j=\b^t_j$ (note that $j=2p_t-2\p-h$ with 
$h\in[1,p_t-\p]$). From the definition we see that $\ta^t_j$ is independent of the choice of $x$, as long as 
$p_x=\p$. For $t\in I$ we set 
$$\nu_t=\sum_{k\in[0,2\p]}n_k|_{k+2p_a-2\p}^a:^t_{2p_t-\p}|
+\sum\Sb r\in I;\\i\in[0,2p_r-2\p-1];\\k\in[0,2\p]\eSb n_k|_{k+i}^r:^t_{2p_t-\p}|\ta^r_i,$$
$$\mu=2/\sqrt{2\nu_a}\qua (\text{if }\nu_a\ne0),\qua\mu=0\qua (\text{if }\nu_a=0).$$
(One can show that $\nu_a$ is nonzero but we will not use this.)

For $t\in I$ and $i,j\in\ZZ$ we set
$$|_i^t:_j^x|=0 \text{ if }-\p\le i-j<2p_t-\p;$$
$$|_i^t:_j^x|=\mu\nu_t \text{ if }i-j=2p_t-\p.$$
For $t\in I$, $i,j\in\ZZ$ such that $i-j=2p_t-\p+s$ $(s\in\ZZ_{>0})$ we define $|_i^t:_j^x|$ by induction on $s$
as follows:
$$\align&|_i^t:_j^x|+\sum_{k\in[1,2\p];k\le s}n_k|_{i-k}^t:_j^x|\\&
=\mu(\sum_{k\in[0,2\p]}n_k|_{i-k}^t:_{j+2p_a-2\p}^a|
+\sum\Sb r\in I;\\h\in[0,2p_r-2\p-1];\\k\in[0,2\p]\eSb n_k|_{i-k}^t:_{j+h}^r|\ta^r_h).\endalign$$
(The right hand side is already known; if $k=s\in[1,2\p]$, the quantity $|_{i-k}^t:_j^x|$ is also known.)
For $t\in I$, $i,j\in\ZZ$ such that $i-j=-\p-s$ ($s\in\ZZ_{>0}$) we define $|_i^t:_j^x|$ by induction on $s$ as
follows:   
$$\align&|_i^t:_j^x|+\sum_{k\in[0,2\p-1];k\ge2\p-s}n_k|_{i+2\p-k}^t:_j^x|\\&
=\mu(\sum_{k\in[0,2\p]}n_k|_{k+2p_a-2\p}^a:_{i-j+2\p}^t|
+\sum\Sb r\in I;\\h\in[0,2p_r-2\p-1];\\k\in[0,2\p]\eSb n_k|_{k+h}^r:_{i-j+2\p}^t|\ta^r_h).\endalign$$
(The right hand side is already known; if $k=2\p-s\in[0,2\p-1]$, the quantity $|_{i+2\p-k}^t:_j^x|$ is also 
known.)
Thus $|_i^y:_j^x|$ is defined for all $y$ such that $p_y>\p$ and $i,j\in\ZZ$ (it is independent of the choice of 
$x$, as long as $p_x=\p$). 

\subhead 2.4 \endsubhead
Assume that $\p\ge1$ and that $p_1>\p$. Then $a,b,I_\p$ are defined (see 2.1). Assume further that $b$ is odd. 
We assume that $y=x$. We write $I$ instead of $I_\p$. 
We set $|^x_i:^x_j|=0$ if $0\le j-i<\p$, $|^x_i:^x_j|=1$ if $j-i=\p$. 

For $i,j\in\ZZ$ such that $j-i=\p+s$ ($s\in\ZZ_{>0}$) we define $|_i^x:_j^x|$ by induction on $s$ as follows:
$$\align&|_i^x:_j^x|+\sum_{k\in[1,2\p];k\le s}n_k|_{i+k}^x:_j^x|\\&
=\mu(\sum_{k\in[0,2\p]}n_k|_{i+k+2p_a-2\p}^a:_j^x|
+\sum\Sb r\in I;\\h\in[0,2p_r-2\p-1];\\k\in[0,2\p]\eSb n_k|_{i+k+h}^r:_j^x|\ta^r_h).\endalign$$
(The right hand side is already known from 2.3; if $k=s\in[1,2\p]$, the quantity $|_{i+k}^x:_j^x|$ is also known.)
Thus $|_i^x:_j^x|$ is defined for $i,j\in\ZZ$, $j\ge i$.

\subhead 2.5 \endsubhead
Assume that $\p\ge1$ and that either $p_1>\p$ and $b$ (see 2.1) is even or that $p_1=\p$. 
We set $|^x_i:^x_j|=0$ if $0\le j-i<\p$, $|^x_i:^x_j|=1$ if $j-i=\p$, $|^x_i:^x_j|=2\p+2$ if $j-i=\p+1$.

For $i,j\in\ZZ$ such that $j-i=\p+s$ ($s\in\ZZ_{\ge2}$) we define $|_i^x:_j^x|$ by induction on $s$ as follows:
$$|^x_i:^x_j|+\sum_{k\in[1,2\p+1];k\le s}(-1)^k\bin{2\p+1}{k}|^x_i:^x_{j-k}|=0.$$ 
We show that, if $i-j-\p=s\ge0$ or if $j-i-\p=s\ge0$, then
$$|_i^x:_j^x|=2(2\p+1)(2\p+2)\do(2\p+s-1)(\p+s)(s!)\i.$$
It is enough to show that
$$\sum_{k\in[0,2\p+1];k\le s}(-1)^k\bin{2\p+1}{k}(2\p+1)...(2\p+s-k-1)(2\p+2s-2k)((s-k)!)\i=0$$
for $s\ge2$ or that
$$\sum_{s\ge0}\sum\Sb k\in[0,2\p+1];\\u\ge0;\\u+k=s\eSb 
(-1)^k\bin{2\p+1}{k}(2\p+1)...(2\p+u-1)(2\p+2u)(u!)\i T^s=1+T$$
or that 
$$\sum_{u\ge0}(2\p+1)...(2\p+u-1)(2\p+2u)(u!)\i T^u=(1+T)(1-T)^{-2\p-1}.$$
More generally we show that
$$\sum_{u\ge0}(M+1)...(M+u-1)(M+2u-1)(u!)\i T^u=(1+T)(1-T)^{-M}$$
for $M\ge2$. The right hand side is equal to 
$$\align&(1-T)^{-M+1}+2T(1-T)^{-M}\\&=\sum_{m\ge0}\bin{M+m-2}{m}T^m+2\sum_{m\ge0}\bin{M+m-1}{m}T^{m+1}\\&=
\sum_{m\ge0}(\bin{M+m-2}{m}+2\bin{M+m-2}{m-1})T^m\\&=\sum_{m\ge0}M(M+1)\do(M+m-2)(M-1+2m)(m!)\i T^m\endalign$$
as desired. Thus $|^x_i:^x_j|$ is defined for all $i,j\in\ZZ$ such that $j\ge i$.

\subhead 2.6\endsubhead
Assume that $\p\ge1$ and that $\p=p_y=p_x$, $y<x$. We set $|^y_i:^x_j|=0$ if $-\p\le i-j<\p$. 

For $i,j\in\ZZ$ such that $j-i=\p+s$ ($s\in\ZZ_{>0}$) we define $|^y_i:^x_j|$ by induction on $s$ as follows:
$$|^y_i:^x_j|+\sum_{k\in[1,2\p];k<s}n_k|^y_{i+k}:^x_j|=\sum_{k\in[0,2\p]}n_k|^x_{i+k}:^x_j|.$$
(The right hand side is known from 2.4, 2.5.)

For $i,j\in\ZZ$ such that $i-j=\p+s$ ($s\in\ZZ_{\ge0}$) we define $|^y_i:^x_j|$  by induction on $s$ as follows:
$$|^y_i:_j^x|+\sum_{k\in[0,2\p-1];k\ge2\p-s}n_k|^y_{i+k-2\p}:_j^x|=\sum_{k\in[0,2q]}n_k|^x_{i+k-2\p}:^x_j|.$$
(The right hand side is known from 2.4, 2.5.)

Under the assumption that either $p_1>\p$ and $b$ is even or that $p_1=\p$ one can show that
$$|_i^y:_j^x|=2\bin{2\p+s}{s} \text{ if } i-j-\p=s\ge0 \text{ or if }j-i-\p-1=s\ge0.$$
This fact will not be used here.

Thus $|^y_i:^x_j|$ is defined for all $i,j\in\ZZ$.

\subhead 2.7\endsubhead
Assume that $\k=1$ and $\p=1/2$. Assume further that $\s$ is even. We have $x=\s+1$. Assume that $y\le\s$. We set
$|^y_i:^x_j|=0$ for all $i,j$. We set $|^x_i:^x_j|=2$ for all $i,j$.

\subhead 2.8\endsubhead
Assume that $\k=1$ and $\p=1/2$. Then $a,b,I_\p$ are defined (see 2.1). We have $b=\s$. Assume further that $\s$ 
is odd. We define $c^r_i$ for $r\in[a,\s],i\in[0,2p_r-1]$ as follows.
If the symmetric matrix (whose entries are already defined)
$$(|^r_i:^{r'}_{i'}|)_{r,r'\in[a,\s],i\in[0,2p_r-1],i'\in[0,2p_{r'}-1]}$$ 
is nonsingular then $c^r_i$ are uniquely defined by the system of linear equations
$$|^a_{2p_a}:^{r'}_{i'}|=\sum_{r\in[a,\s],i\in[0,2p_r-1]}|^r_i:^{r'}_{i'}|c^r_i$$
whose coefficients are already defined. If the symmetric matrix above is singular then we set $c^r_i=0$ for
all $r,i$. (One can show that the last possibility does not occur; this will not be used here.) We set
$$c^x_0=\sqrt{\nu/2}$$
where
$$\nu=-\sum_{r,r'\in[a,\s];i\in[0,2p_r-1];i'\in[0,2p_{r'}-1]}c^r_ic^{r'}_{i'}|^r_i:^{r'}_{i'}|$$
is already defined. For any $r'\in[a,\s],i'\in[0,2p_{r'}-1],h\in\ZZ$ we set
$$|^{r'}_{i'}:^x_h|=
(c^x_0)\i(|^a_{2p_a+h}:^{r'}_{i'}|-\sum_{r\in[a,\s],i\in[0,2p_r-1]}c^r_i|^r_{i+h}:^{r'}_{i'}|)$$
(the right hand side is already defined).
For $h,h'\in\ZZ$ we set
$$\align&|^x_h:^x_{h'}|=(c^x_0)^{-2}(|^a_{2p_a+h}:^a_{2p_a+h'}|
-\sum\Sb r\in[a,\s];\\i\in[0,2p_r-1]\eSb c^r_i(|^r_{i+h}:^a_{2p_a+h'}|+|^r_{i+h'}:^a_{2p_a+h}|\\&
+\sum_{r,r'\in[a,\s];i\in[0,2p_r-1];i'\in[0,2p_{r'}-1]}c^r_ic^{r'}_{i'}|^r_{i+h}:^{r'}_{i'+h'}|)\endalign$$
(the right hand side is already defined).

\subhead 2.9\endsubhead
This completes the inductive definition of $|^y_i:^x_j|$. From the definitions we see that
$|^y_i:^x_j|=|^y_{i'}:^x_{j'}|$ if $i-j=i'-j'$.

\subhead 2.10\endsubhead
Assume for example that $\s=2$, $k=p_1>p_2=1$. It is likely that $|^1_{2k}:2^1|^2=(-1)^{k-1}2^{2k}$. (This is 
true at least if $k\in\{2,3,4\}$.) In particular, if $\kk=\CC$, $|^y_i:^x_j|$ is not necesarily a real number.

\subhead 2.11\endsubhead
In the remainder of this section we fix $V,Q,(,),\nn,\k,Is(V)$ as in 1.1; we assume that $Q\ne0$, $p\ne2$.
We shall assume that $p_1+p_2+\do+\p_\s=(\nn-\k)/2$; if $\k=0$ we assume also that $\k_\s=0$. When $\k=1$ we set 
$p_{\s+1}=1/2$. 

\subhead 2.12\endsubhead
We fix a unipotent element $g\in Is(V)$. Let $N=g-1:V@>>>V$. We shall assume that $\cm(N,V)$ consists of 
$2p_1+\ps(1)\ge2p_2+\ps(2)\ge\do\ge2p_\s+\ps(\s)$ (and $1$, if $\k=1$ and $\k_\s=0$). 

\proclaim{Proposition 2.13} Let $w^t_i$ ($t\in[1,\s+\k],i\in\ZZ$) be a $(g,p_*)$-adapted collection of vectors in
$V$ (see 1.2). Then there exists $\e:[1,\s+\k]@>>>\{1,-1\}$, $t\m\e_t$ such that 
$(\e_tw^t_i,\e_rw^r_j)=|^t_i:^r_j|$ (see 2.9) for all $t,r\in[1,\s+\k],i,j\in\ZZ$.
\endproclaim
The proof is given in 2.17-2.26.

\subhead 2.14\endsubhead
Note that

(a) {\it the collection $\{N^iw^r_0; r\in[1,\s+\k],i\in[0,2p_r-1]\}$ is a basis of $V$.}
\nl
This follows from 1.3(a) since the collection in (a) is related to the collection in 1.3(a) by an upper triangular
matrix with $1$ on diagonal. Similarly, the following statement can be deduced from 1.3(b). Let 
$e,f\in[1,\s+\k]$, $e\le f$ be such that the subspace of $V$ spanned by $\{N^iw^x_0;x\in[e,f],i\in[0,2p_x-1]\}$ 
is $N$-stable; then

(b) {\it the symmetric matrix $(N^iw^x_0,N^jw^y_0)_{x,y\in[e,f],i\in[0,2p_x-1],j\in[0,2p_y-1]}$ is nonsingular.}

\subhead 2.15\endsubhead
Let $r\in[1,\s]$ be such that $\ps(r)=-1$. Then $\sum_{t\in[1,r]}\ps(t)=0$ hence
$$\dim(N^{2p_r-1}V)=\sum_{t\in[1,r]}(2p_t+\ps(t)-2p_r+1)=\sum_{t\in[1,r]}(2p_t-2p_r+1).$$
Let $\cw=\cw_{1,r}$. Let $\cw'=\cw_{t+1,\s+\k}$ (convention: if $\k=0,t=\s$ then $\cw'=0$). We show:

(a) $g\cw=\cw$; $g\cw'=\cw'$; $\cw'=\cw^\pe$; 

(b) {\it $\cm(N,\cw)$ consists of $2p_1+\ps(1)\ge2p_2+\ps(2)\ge\do\ge2p_r+\ps(r)$; $\cm(N,\cw')$ consists of
$2p_{r+1}+\ps(r+1)\ge2p_{r+2}+\ps(r+2)\ge\do\ge2p_\s+\ps(\s)$ (and $1$, if $\k=1$ and $\k_\s=0$).}
\nl
(a) is a special case of 2.4(a) (with $k=2p_r-1$, $d=r$). Note that $2p_r>2p_r-1\ge 2p_{r+1}$.

We prove (b). From 2.4(a) we see that $N:\cw@>>>\cw$ has exactly $r$ Jordan blocks. These are some of the Jordan 
blocks of $N:V@>>>V$ (recall that $V=\cw\op\cw'$ and $\cw,\cw'$ are $N$-stable). Hence $\cm(N,\cw)$ is given by 
$r$ terms of the sequence 
$2p_1+\ps(1)\ge2p_2+\ps(2)\ge\do\ge2p_\s+\ps(\s)$ (and $1$, if $\k=1$ and $\k_\s=0$). Now the sum of numbers in
$\cm(N,\cw)$ is equal to $\dim\cw$ which is equal to the sum of the first $r$ terms of the sequence above. Hence 
these numbers must be given by the first $r$ terms of our sequence. This proves the first statement of (b). Using
again that $V=\cw\op\cw'$ and $\cw,\cw'$ are $N$-stable, we see that the second statement of (b) follows from the 
first statement of (b).

\subhead 2.16\endsubhead
Let $c\in[1,\s]$ be such that $\ps(c)=1$. (Thus $c$ is odd.) Let $\p$ be such that $\p=p_x$ for some $x\in[1,\s]$
and $p_c\ge\p$. Let $I=\{y\in[c,\s];p_y>\p\}$. Note that $I=\em$ if $p_c=\p$ and $I$ is of the form $[c,d]$ 
with $c\le d\le\s$ if $p_c>\p$; in this last case we assume that $\ps(t)=0$ for any $t\in I-\{c\}$. If $p_c>\p$
then $d$ is odd. (If $d<\s$ then $p_d>p_{d+1}$, $\ps(d)=0$, hence $d$ is odd. If $d=\s$ is even then 
$\ps(d)=-1$, contrary to our assumption.) Let $\cz$ be the subspace of $V$ spanned by 
$\{N^iw^y_0;y\in I,i\in[2\p,2p_y-1]\}$ and by $N^{2p_c}w^c_0$. Equivalently, $\cz$ is the subspace of $V$ spanned
by $\{N^{2\p}w^y_i;y\in I,i\in[0,2p_y-2\p-1]\}$ and by $N^{2\p}w^c_{2p_c-2\p}$. We show:

(a) $N^{2\p}w^x_0\in\cz$ for any $x\in[c,\s]$ such that $p_x=\p$;

(b) $\dim\cz=1+\sum_{y\in I}(2p_y-2\p)$.
\nl
Let $\cw'$ be the subspace of $V$ spanned by $\{w^z_i;z\in[c,\s+\k],i\in[0,2p_z-1]\}$ or equivalently by 
$\{N^iw^z_0;z\in[c,\s+\k],i\in[0,2p_z-1]\}$. (This agrees with the notation of 2.15 if $c>1$ and $r=c-1$; if $c=1$
we have $\cw'=V$.) By 2.15, $\cw'$ is $N$-stable and $\cm(N,\cw')$ consists of
$2p_c+\ps(c)\ge2p_{c+1}+\ps(c+1)\ge\do\ge2p_\s+\ps(\s)$ (and $1$, if $\k=1$ and $\k_\s=0$). In particular,
$\dim N^{2p_c}\cw'=1$.

We prove (a) and (b) by induction on $p_c-\p$. 
Assume first that $p_c=\p$ so that $I=\em$. Assume that (b) is false that is, $N^{2p_c}w^c_0=0$. Then the 
subspace $X$ of $V$ spanned by $N^iw^c_0$ ($i\in[0,2p_c-1])$ is $N$-stable. By 2.14(b), $(,)|_X$ is 
nonsingular. The vectors $\{N^iw^c_0;i\in[0,2p_c-1]\}$ are linearly independent by 2.14(a). Hence $N+1:X@>>>X$
is a unipotent isometry of $X$ with a single Jordan block and $X$ has even dimension $>0$; this is impossible. 
Thus (b) holds. As we have just seen we have $N^{2p_c}w^c_0\ne0$. Since $\dim N^{2p_c}\cw'=1$ we see that
$N^{2p_c}\cw'$ is spanned by $N^{2p_c}w^c_0$. Since $N^{2p_c}\cw'=1$ is an $N$-stable line we see that

(c) $N^{2p_c+1}\cw'=0$.
\nl
Now let $x\in[c,\s]$ be such that $p_x=\p$. Then $N^{2p_x}w^x_0\in N^{2p_c}\cw'$ hence $N^{2p_x}w^x_0$ is a 
multiple of $N^{2p_c}w^t_0$ and thus is in $\cz$. We see that (a) holds. 

In the rest of the proof we assume that $p_c>\p$ and that (a) and (b) hold for any $\p'$ such that $\p'=p_r$ for 
some $r\in[1,\s]$ and $p_c\ge\p'>\p$. 
Let $\cz'$ be the subspace of $V$ spanned by $\{N^iw^y_0;y\in I,i\in[2\p,2p_y-1]\}$. These vectors are linearly
independent by 2.14(a). Hence $\dim\cz'=\sum_{y\in I}(2p_y-2\p)$. We see that (b) is equivalent to the equality 
$\dim\cz=1+\dim\cz'$. Assume that this equality is not true. Then $N^{2p_c}w^c_0\in\cz'$. We show that $\cz'$ is 
$N$-stable. If $y\in I$, we have $p_c\ge p_y>\p$ and the induction hypothesis shows that $N^{2p_y}w_0^y$ is a 
linear combination of $\{N^iw^{y'}_0;y'\in[c,y],i\in[2p_y,2p_y-1]\}$ (which are in $\cz'$) and of $N^{2p_c}w^c_0$
(which is also in $\cz'$). Thus, $N^{2p_y}w_0^y\in\cz'$. We see that $N$ maps each of the basis elements of 
$\cz'$ to another element of that basis or to $N^{2p_y}w_0^y\in\cz'$, $(y\in I)$. Thus, $N\cz'\sub\cz'$. Let 
$\cz''$ be 
the subspace of $V$ spanned by $\{N^iw^y_0;y\in I,i\in[0,2p_y-1]\}$. We have $\cz'\sub\cz''$. We show that 
$N\cz''\sub\cz''$. It is enough to show that for any $y\in I$ we have $N^{2p_y}w_0^y\in\cz''$. This follows from 
$N^{2p_y}w_0^y\in\cz'$ and $\cz'\sub\cz''$. By 2.14(b), $(,)|_{\cz''}$ is nonsingular. Hence we have 
$\cw'=\cz''\op\ti\cz''$ where $\ti\cz''=\{x\in\cw';(x,\cz'')=0\}$. Moreover, since $\cw',\cz''$ are $g$-stable we
see that $\ti\cz''$ is also $g$-stable hence $N$-stable. Hence  (setting $N_0=N|_{\cz''}$) we see that 
$\cm(N_0,\cz'')$ is contained in $\cm(N,\cw')$. Thus,

(d) {\it $\cm(N_0,\cz'')$ is contained in the multiset which consists of the numbers\lb
$2p_c+\ps(c)\ge2p_{c+1}+\ps(c+1)\ge\do\ge2p_\s+\ps(\s)$ (and $1$, if $\k=1$ and $\k_\s=0$).}
\nl
By the first part of the proof, the vectors $\{N^iw^c_0;i\in[0,2p_c]\}$ of $\cz''$ are linearly independent.
Hence $N_0$ has at least one Jordan block of size $\ge2p_c+1$ and using (d) it has exactly one Jordan block of 
size $2p_c+1$. Applying \cite{\WE, 3.1} to $N_0$ we see that $N_0$ has at most $|I|=d-c+1$ Jordan blocks. If 
$c=d$ then it follows that $N_0$ is a single Jordan block and we must have $\dim\cz''=2p_c+1$; but in this case 
we have $\dim\cz''=2p_c$, contradiction. Now assume that $c<d$. Then $\cm':=\cm(N_0,\cz'')-\{2p_c+1\}$ is a part 
of the list $2p_{c+1}\ge\do\ge2p_d\ge\do$. Hence if $S$ is the sum of the 
numbers in $\cm'$ (that is, $S=\dim\cz''-(2p_c+1)=(2p_c+2p_{c+1}+\do+2p_d)-(2p_c+1)=2p_{c+1}+\do+2p_d-1$), then
$S\le2p_{c+1}+\do+2p_d$. This implies that the numbers in $\cm'$ are obtained from the list 
$2p_{c+1}\ge\do\ge2p_d$ by decreasing exactly one number in the list by $1$ and leaving the other numbers
unchanged. It follows that $N_0$ has exactly $d-c+1$ Jordan blocks of which two have odd size and the remaining 
$d-c-1$ are of even size. But $c$ and $d$ are odd hence $d-c-1$ is odd. Thus $N_0$ has an odd number of Jordan 
blocks of even size. This is not possible since $N_0+1$ is a unipotent isometry of the nonsingular form $(,)$ on
the even dimensional space $\cz''$. This contradiction proves (b).

Now let $x\in[c,\s]$ be such that $p_x=\p$. From the knowledge of $\cm(N,\cw')$ we see that 
$\dim N^{2\p}\cw'=1+\sum_{y\in I}(2p_y-2\p)=1+\dim\cz'$. Using (b) we deduce that $\dim N^{2\p}\cw'=\dim\cz$. 
From the definitions we have $\cz\sub N^{2\p}\cw'$. It follows that $\cz=N^{2\p}\cw'$. Clearly, we have 
$N^{2p_x}w^x_0\in N^{2\p}\cw'$. Hence $N^{2p_x}w^x_0\in\cz$ and (a) is proved.

This completes the inductive proof of (a) and (b).

\subhead 2.17\endsubhead
Since $(w^t_i,w^r_j)=(w^r_j,w^t_i)$, to prove 2.13, it is enough to prove that 

$(\e_yw^y_i,\e_xw^x_j)=|^y_i:^x_j|$
\nl
under the following assumptions (which will be in force until the end of 2.26):

$\p$ is a fixed number equal to one of $p_1,p_2,\do,p_{\s+\k}$;

$p_x=\p$, $y\le x$ (and $i\le j$ if $y=x$);

$(\e_tw^t_i,\e_rw^r_j)=|^t_i:^r_j|$ whenever $p_t>\p,p_r>\p$.
\nl
Here $\e_t\in\{1,-1\}$ are already defined for all $t$ such that $p_t>\p$ and $\e_t\in\{1,-1\}$ is to be defined 
for all $t$ such that $p_t=\p$. For all $t$ such that $p_t>\p$ we replace $w^t_i$ by $\e_tw^t_i$ and we see that 
we can assume that $\e_t=1$ for such $t$. Thus we have $(w^t_i,w^r_j)=|^t_i:^r_j|$ whenever $p_t>\p,p_r>\p$.

In the case where $p_1>\p$ we define integers $a,b\in[1,\s]$ as in 2.1; we set $I_\p=[a,b]$.
For any $k\in[0,2\p]$ let $n_k$ be as in 2.1.

\subhead 2.18\endsubhead
Assume that $p_y>\p$ and that for some even $r\in[y,x-1]$ we have $p_r>p_{r+1}$. We then have $\ps(r)=-1$. Define 
$\cw,\cw'$ in terms of $r$ as in 2.15. We have $w^y_i\in\cw,w^x_j\in\cw'$, $(\cw,\cw')=0$ (see 2.15(a)). Hence we
have $(w^y_i,w^x_j)=0$. Thus 2.13 holds in this case.

\subhead 2.19\endsubhead
Assume that $p_y>\p\ge1$ and that for any even $r\in[y,x-1]$ we have $p_r=p_{r+1}$. Note that $p_1>\p$ hence 
$I_\p=[a,b]$ is defined as in 2.17. Now the assumptions of 2.16 are satisfied with $c=a$, $I=I_\p$, $d=b$. We 
have $y\in I$. Using 2.16(a) we see that there exist 
$A^r_i\in\kk$ $(r\in I;i\in[0,2p_r-2\p-1])$ and $C\in\kk$ such that
$$N^{2\p}w^x_0=CN^{2\p}w_{2p_a-2\p}^a+\sum_{r\in I;i\in[0,2p_r-2\p-1]}A^r_iN^{2\p}w^r_i.$$
Thus we have
$$\sum_{k\in[0,2\p]}n_kw_k^x=C\sum_{k\in[0,2\p]}n_kw_{k+2p_a-2\p}^a
+\sum_{r\in I;i\in[0,2p_r-2\p-1];k\in[0,2\p]}A^r_in_kw_{k+i}^r.\tag a$$
Applying $(,w^x_\p)$ to (a), we obtain
$$\align&\sum_{k\in[0,2\p]}n_k(w_k^x,w^x_\p)=C\sum_{k\in[0,2\p]}n_k(w_{k+2p_a-2\p}^a,w^x_\p)\\&
+\sum_{r\in I;i\in[0,2p_r-2\p-1];k\in[0,2\p]}A^r_in_k(w_{k+i}^r,w^x_\p).\endalign$$
In the first sum we have $(w_k^x,w^x_\p)=0$ except when $k=0$ or $k=2\p$. In the second sum we have
$k+2p_a-2\p\in[0,2p_a-1]$ hence $(w_{k+2p_a-2\p}^a,w^x_\p)=0$ except when $k=2\p$. In the third sum we have 
$k+i\in[0,2p_r-1]$ hence each term of the sum is zero. Thus we have

{\it $2=C(w_{2p_a}^a,w^x_p)$; in particular, $C\ne0$.}
\nl
For $r\in I,i\in[0,p_r-\p-1]$ we set $\tA^r_i=C\i A^r_i$; for $r\in I,i\in[p_r-\p,2p_r-2\p-1]$ we set 
$\tB^r_i=C\i A^r_i$; then we have
$$\align&\sum_{k\in[0,2\p]}C\i n_kw_k^x=\sum_{k\in[0,2\p]}n_kw_{k+2p_a-2\p}^a
+\sum\Sb r\in I;\\i\in[0,p_r-\p-1]\eSb\tA^r_i\sum_{k\in[0,2\p]}n_kw_{k+i}^r\\&
+\sum_{r\in I;j\in[1,p_r-\p]}\tB^r_{2p_r-2\p-j}\sum_{k\in[0,2\p]}n_kw_{k+2p_r-2\p-j}^r.\tag b\endalign$$
Apply $(,w^r_u)$ to (b) where $r\in I$, $u\in[\p,2p_r-\p]$ and note that 
$(w_u^r,w_k^x)=(w_{u-k+\p}^r,w_\p^x)$ equals $0$ except when $u=2p_r-\p,k=0$ when it is $(w_{2p_r}^r,w_\p^x)$.
(For $k\in[0,2\p]$ we have $0\le2\p-k\le u-k+\p\le u+\p\le2p_r$.)
We also substitute $(w_i^t,w_j^{t'})=|_i^t:_j^{t'}|$ for $t,t'\in I$. We obtain
$$\align&\d_{u,2p_r-\p}C\i(w^r_{2p_r},w_\p^x)=\sum_{k\in[0,2\p]}n_k|_{k+2p_a-2\p}^a:^r_u|\\&
+\sum_{r'\in I;i\in[0,p_{r'}-\p-1]}\tA^{r'}_i\sum_{k\in[0,2\p]}n_k|_{k+i}^{r'}:^r_u|\\&
+\sum_{r'\in I;j\in[1,p_{r'}-\p]}\tB^{r'}_{2p_{r'}-2\p-j}\sum_{k\in[0,2\p]}n_k|_{k+2p_{r'}-2\p-j}^{r'}:^r_u|.
\tag b1\endalign$$
We take $u=p_r+h$ ($h\in[0,p_r-\p-1]$) or $u=p_r-h$ ($h\in[1,p_r-\p]$). We obtain
$$\align&0=\sum_{k\in[0,2\p]}n_k|_{k+2p_a-2\p}^a:^r_{p_r+h}|
+\sum'\Sb r'\in I;\\i\in[0,p_{r'}-\p-1];\\k\in[0,2\p]\eSb\tA^{r'}_in_k|_{k+i}^{r'}:^r_{p_r+h}|\\&
+\sum''_{r'\in I;j\in[1,p_{r'}-\p];k\in[0,2\p]}\tB^{r'}_{2p_{r'}-2\p-j}n_k|_{k+2p_{r'}-2\p-j}^{r'}:^r_{p_r+h}|
\tag c\endalign$$
for $r\in I,h\in[0,p_r-\p-1]$ and
$$\align&0=\sum_{k\in[0,2\p]}n_k|_{k+2p_a-2\p}^a:^r_{p_r-h}|
+\sum'\Sb r'\in I;\\i\in[0,p_{r'}-\p-1];\\k\in[0,2\p]\eSb\tA^{r'}_in_k|_{k+i}^{r'}:^r_{p_r-h}|\\&
+\sum''_{r'\in I;j\in[1,p_{r'}-\p];k\in[0,2\p]}\tB^{r'}_{2p_{r'}-2\p-j}n_k|_{k+2p_{r'}-2\p-j}^{r'}:w^r_{p_r-h}|
\tag d\endalign$$  
for $r\in I,h\in[1,p_r-\p]$.

In $\sum'$ for (c) we have assuming $i>h$:

if $r'\le r$ then $1\le k+i-h\le2\p+p_{r'}-\p-1\le2p_{r'}-1$ hence $|_{k+i}^{r'}:^r_{p_r+h}|=0$;

if $r'\ge r$ then $1\le p_r-2\p-p_{r'}+\p+1+p_{r'}\le p_r+h-k-i+p_{r'}\le p_r-1+p_{r'}\le2p_r-1$ hence 
$|_{k+i}^{r'}:^r_{p_r+h}|=0$.

In $\sum'$ for (c) we have assuming $i=h$:

if $r'\le r$ then $0\le k+i-h\le2\p+p_{r'}-\p-1\le2p_{r'}-1$ hence $|_{k+i}^{r'}:^r_{p_r+h}|$ is $0$ unless 
$r'=r,k=0$ when it is $1$;

if $r'>r$ then $0\le p_r-2\p+p_{r'}\le p_r+h-k-i+p_{r'}\le p_r+p_{r'}\le2p_r$ hence $|_{k+i}^{r'}:^r_{p_r+h}|$ is
$0$ unless $k=0,p_r=p_{r'}$ when it equals $|_0^{r'}:^r_{p_r}|$.

In $\sum''$ for (c) we have assuming $j\ge h$ (so that $j+h\le2j-\d_{j,p_{r'}-\p}$):

if $r'<r$ then $0\le2p_{r'}-2\p-2j\le k+2p_{r'}-2\p-j-h\le2\p+2p_{r'}-2\p-1=2p_{r'}-1$ hence 
$|_{k+2p_{r'}-2\p-j}^{r'}:^r_{p_r+h}|=0$;

if $r'\ge r$ then 
$1\le p_r-2\p-p_{r'}+2\p+1\le p_r+h-k-p_{r'}+2\p+j\le p_r+2j-\d_{j,p_{r'}-\p}-p_{r'}+2\p\le 
p_r-p_{r'}+2\p+2p_{r'}-2\p-1=p_r+p_{r'}-1\le2p_r-1$ hence $|_{k+2p_{r'}-2\p-j}^{r'}:^r_{p_r+h}|=0$.

In $\sum'$ for (d) we have assuming $i\ge h-1$ (hence $i+h\le2i+1$):

if $r'\le r$ then $1\le k+i+h\le2\p+2i+1\le2\p+2p_{r'}-2\p-2+1\le2p_{r'}-1$ hence $|_{k+i}^{r'}:^r_{p_r-h}|=0$;

if $r'\ge r$ then
$1\le p_r-p_{r'}+1=p_r-2\p+p_{r'}-2p_{r'}+2\p+2-1\le p_r-2\p+p_{r'}-2i-1\le p_r-h-k-i+p_{r'}\le p_r+p_{r'}-1
\le2p_r-1$ hence $|_{k+i}^{r'}:^r_{p_r-h}|=0$.

In $\sum''$ for (d) we have assuming $j>h$:

if $r'\le r$ then $1\le 2p_{r'}-2\p-p_{r'}+\p+1\le k+2p_{r'}-2\p-j+h\le2\p+2p_{r'}-2\p-1=2p_{r'}-1$ hence 
$|_{k+2p_{r'}-2\p-j}^{r'}:^r_{p_r-h}|=0$;

if $r'\ge r$ then $1\le p_r+1-p_{r'}\le p_r-h-k-p_{r'}+2\p+j\le p_r-1-p_{r'}+2\p+p_{r'}-\p=p_r-1+\p\le2p_r-1$
hence $|_{k+2p_{r'}-2\p-j}^{r'}:^r_{p_r-h}|=0$.

In $\sum''$ for (d) we have assuming $j=h$:

if $r'<r$ then $0\le 2p_{r'}-2\p\le k+2p_{r'}-2\p-j+h\le2\p+2p_{r'}-2\p=2p_{r'}$ hence 
$|_{k+2p_{r'}-2\p-j}^{r'}:^r_{p_r-h}|$ is $0$ unless $k=2\p$ when it equals $|_{2p_{r'}}^{r'}:^r_{p_r}|$;

if $r'\ge r$ then $0\le p_r-p_{r'}\le p_r-h-k-p_{r'}+2\p+j\le p_r-p_{r'}+2\p\le2p_r-1$ hence 
$|_{k+2p_{r'}-2\p-j}^{r'}:^r_{p_r-h}|$ is $0$ unless $r'=r,k=2\p$ when it equals $1$.
\nl
Thus (c),(d) can be rewritten as follows (we also substitute $(w_i^t,w_j^{t'})=|_i^t:_j^{t'}|$ for $t,t'\in I$):
$$\align&\tA^r_h+\sum'_{r'\in I;r'>r;p_r=p_r'}\tA^{r'}_h|_0^{r'}:^r_{p_r}|
=-\sum_{k\in[0,2\p]}n_k|_{k+2p_a-2\p}^a:^r_{p_r+h}|\\&
-\sum'_{r'\in I;i\in[0,p_{r'}-\p-1];i<h;k\in[0,2\p]}\tA^{r'}_in_k|_{k+i}^{r'}:^r_{p_r+h}|\\&
-\sum''_{r'\in I;j\in[1,p_{r'}-\p];j<h;k\in[0,2\p]}\tB^{r'}_{2p_{r'}-2\p-j}n_k|_{k+2p_{r'}-2\p-j}^{r'}:^r_{p_r+h}|
\tag c1\endalign$$
for $r\in I,h\in[0,p_r-\p-1]$ and
$$\align&\tB^r_{2p_r-2\p-h}+\sum''_{r'\in I;r'<r}\tB^{r'}_{2p_{r'}-2\p-h}|_{2p_{r'}}^{r'}:^r_{p_r}|
=-\sum_{k\in[0,2\p]}n_k|_{k+2p_a-2\p}^a:^r_{p_r-h}|\\&
-\sum'_{r'\in I;i\in[0,p_{r'}-\p-1];i<h-1;k\in[0,2\p]}\tA^{r'}_in_k|_{k+i}^{r'}:^r_{p_r-h}|\\&
-\sum''_{r'\in I;j\in[1,p_{r'}-\p];j<h;k\in[0,2\p]}\tB^{r'}_{2p_{r'}-2\p-j}n_k|_{k+2p_{r'}-2\p-j}^{r'}:^r_{p_r-h}|
\tag d1\endalign$$  
for $r\in I,h\in[1,p_r-\p]$.
Note that (c1),(d1) can be viewed as inductive formulas for $\tA^r_h,\tB^r_{2p_r-2\p-h}$ which are identical to 
the inductive formulas 2.3(i),(ii). It follows that $\tA^r_h=\a^r_h$ for $r\in I,h\in[0,p_r-\p-1]$ and 
$\tB^r_{2p_r-2\p-h}=\b^r_{2p_r-2\p-h}$ for $r\in I,h\in[1,p_r-\p]$. We define $\ta^t_j$ for 
$t\in I,j\in[0,2p_t-2\p-1]$ as in 2.3. Now (b1) with $u=2p_r-\p$ becomes
$$\align&C\i(w^r_{2p_r},w_\p^x)=\sum_{k\in[0,2\p]}n_k|_{k+2p_a-2\p}^a:^r_{2p_r-\p}|\\&
+\sum_{r'\in I;i\in[0,2p_{r'}-2\p-1]}\ta^{r'}_i\sum_{k\in[0,2\p]}n_k|_{k+i}^{r'}:^r_{2p_r-\p}|.\endalign$$
Equivalently,
$$C\i(w^r_{2p_r},w_\p^x)=\nu_r$$
with $\nu_r$ as in 2.3. Combining this with the earlier identity $2=C(w_{2p_a}^a,w^x_\p)$ we see that 
$2\nu_r=(w^r_{2p_r},w_\p^x)(w_{2p_a}^a,w^x_\p)$. Taking $r=a$ we see that $2\nu_a=(w_{2p_a}^a,w^x_\p)^2$. Hence 
$(w_{2p_a}^a,w^x_\p)=\e_x\sqrt{2\nu_a}$ where $\e_x\in\{1,-1\}$. Replacing $w^x_i$ by $\e_xw^x_i$ for $i\in\ZZ$ 
we see that we can assume that $\e_x=1$ so that $(w_{2p_a}^a,w^x_\p)=\sqrt{2\nu_a}$. From 
$2=C(w_{2p_a}^a,w^x_\p)=C\sqrt{2\nu_a}$ we see that $\nu_a\ne0$ and $C=2/\sqrt{2\nu_a}$. Thus with the notation 
of 2.3 we have $C=\mu$. We deduce $(w^r_{2p_r},w_\p^x)=\mu\nu_r$ hence $(w^r_{2p_r},w_\p^x)=|^r_{2p_r}:_\p^x|$ 
for $r\in I$. 

We can write (b) as follows 
$$\sum_{k\in[0,2\p]}\mu\i n_kw_k^x=\sum_{k\in[0,2\p]}n_kw_{k+2p_a-2\p}^a
+\sum\Sb r\in I;\\i\in[0,2p_r-2\p-1]\eSb\ta^r_i\sum_{k\in[0,2\p]}n_kw_{k+i}^r.\tag e$$
Applying $(,w^r_{2p_r-\p+s})$ (where $s\in\ZZ_{>0})$ to (e) and using the induction hypothesis we obtain
$$\align&\mu\i(w^r_{2p_r-\p+s},w_0^x)+\sum_{k\in[1,2\p];k\le s}\mu\i n_k|^r_{2p_r-\p+s}:_k^x|\\&
=\sum_{k\in[0,2\p]}n_k|_{k+2p_a-2\p}^a:^r_{2p_r-\p+s}|\\&
+\sum_{r'\in I;i\in[0,2p_{r'}-2\p-1];k\in[0,2\p]}\ta^{r'}_in_k|_{k+i}^{r'}:^r_{2p_r-\p+s}|.\endalign$$
(We use that $(w^r_{2p_r-\p+s},w_k^x)=0$ if $k>s$; indeed, $0\le2p_r-2\p\le2p_r+s-k\le2p_r-1$.) Comparing this 
with 2.3 we obtain
$$\align&(w^r_{2p_r-\p+s},w_0^x)+\sum_{k\in[1,2\p];k\le s}n_k(w^r_{2p_r-\p+s},w_k^x)\\&=
|^r_{2p_r-\p+s}:_0^x|+\sum_{k\in[1,2\p];k\le s}n_k|^r_{2p_r-\p+s}:_k^x|\endalign$$
for $s>0$ which implies by induction on $s$ that $(w^r_{2p_r-\p+s},w_0^x)=|^r_{2p_r-\p+s}:_0^x|$.

Applying $(,w^r_{\p-s})$ (where $s\in\ZZ_{>0})$ to (e) and using the induction hypothesis we obtain
$$\align&\mu\i(w^r_{\p-s},w_{2\p}^x)+\sum_{k\in[0,2\p-1];k\ge2\p-s}\mu\i n_k(w^r_{\p-s},w_k^x)\\&=
\sum_{k\in[0,2\p]}n_k|_{k+2p_a-2\p}^a:^r_{\p-s}|
+\sum_{r'\in I;i\in[0,2p_{r'}-2\p-1]}\ta^{r'}_i\sum_{k\in[0,2\p]}n_k|_{k+i}^{r'}:^r_{\p-s}|.\endalign$$
(We use that $(w^r_{\p-s},w_k^x)=0$ if $k<2\p-s$; indeed $1\le2\p-s-k\le 2p_r-1$.)
Comparing this with 2.3 we obtain
$$\align&(w^r_{\p-s},w_{2\p}^x)+\sum_{k\in[0,2\p-1];k\ge2\p-s}n_k(w^r_{\p-s},w_k^x)\\&=
|^r_{\p-s}:_{2\p}^x|+\sum_{k\in[0,2\p-1];k\ge2\p-s}n_k|^r_{\p-s}:_k^x|\endalign$$
for $s>0$ which implies by induction on $s$ that
$$(w^r_{\p-s},w_{2\p}^x)=|^r_{\p-s}:_{2\p}^x|.$$
We see that 
$$(w^r_i,w_j^x)=|^r_i:_j^x|\tag f$$ 
for any $r\in I$ and any $i,j\in\ZZ$.

\subhead 2.20\endsubhead
Assume that $\p\ge1$ and that $p_1>\p$. Then $a,b,I_\p$ are defined (see 2.17). Assume further that $b$ is odd. 
Then 2.19(e) holds. We write $I$ instead of $I_\p$. Applying $(,w^x_{\p+s})$ where $s\in\ZZ_{>0}$ to 2.19(e) we 
obtain
$$\align&\mu\i (w_0^x,w^x_{\p+s})+\sum_{k\in[1,2\p];k\le s}\mu\i n_k(w_k^x,w^x_{\p+s})\\&
=\sum_{k\in[0,2\p]}n_k(w_{k+2p_a-2\p}^a,w^x_{\p+s})
+\sum_{r\in I;i\in[0,2p_r-2\p-1]}\ta^r_i\sum_{k\in[0,2\p]}n_k(w_{k+i}^r,w^x_{\p+s}).\endalign$$
(We use that $(w_k^x,w^x_{\p+s})=0$ for $k>s$; indeed we have $1\le2\p+s-k\le2\p-1$.) We rewrite this using 
2.19(f):
$$\align&\mu\i(w_0^x,w^x_{\p+s})+\sum_{k\in[1,2\p];k\le s}\mu\i n_k(w_k^x,w^x_{\p+s})\\&
=\sum_{k\in[0,2\p]}n_k|_{k+2p_a-2\p}^a:w^x_{\p+s}|
+\sum_{r\in I;i\in[0,2p_r-2\p-1]}\ta^r_i\sum_{k\in[0,2\p]}n_k|_{k+i}^r:^x_{\p+s}|.\endalign$$
Comparing with 2.4 we deduce
$$(w_0^x,w^x_{\p+s})+\sum_{k\in[1,2\p];k\le s}n_k(w_k^x,w^x_{\p+s})=
|_0^x:^x_{\p+s}|+\sum_{k\in[1,2\p];k\le s}n_k|_k^x:^x_{\p+s}|.$$
From this equality we see by induction on $s$ that $(w_0^x,w^x_{\p+s})=|_0^x:^x_{\p+s}|$.
Hence $(w_i^x,w^x_j)=|_i^x:^x_j|$ for any $i,j$ such that $j-i>\p$. By symmetry the same equality holds 
for any $i,j$ such that $i-j>\p$. It also holds for $j-i\in[-\p,\p]$ by the definition of $w^x_i$. Hence it holds
for any $i,j$.

\subhead 2.21\endsubhead
Assume that $\p\ge1$ and that either $p_1>\p$ and $b$ (see 2.17) is even or that $p_1=\p$. Now 2.16(c) is 
applicable with $c=b+1$ (if $p_1>\p$) or with $c=1$ (if $p_1=\p$). From 2.16(c) we see that $N^{2\p+1}w^x_0=0$ 
that is,
$$\sum_{k\in[0,2\p+1]}(-1)^k\bin{2\p+1}{k}w^x_k=0.\tag a$$
Applying $(,w^x_\p)$ to (a) we obtain
$$\sum_{k\in\{0,2\p,2\p+1\}}(-1)^k\bin{2\p+1}{k}(w^x_k,w^x_\p)=0$$ 
that is $1+(2\p+1)-(w^x_{2\p+1},w^x_\p)=0$ so that $(w^x_{\p+1},w^x_0)=2\p+2$. Thus 
$(w^x_{\p+1},w^x_0)=|^x_{\p+1}:^x_0|$, see 2.5.
Aplying $(,w^x_{\p+s})$ (with $s\in\ZZ_{\ge2}$) to (a) we obtain
$$(w^x_0,w^x_{s+\p})+\sum_{k\in[1,2\p+1];k\le s}(-1)^k\bin{2\p+1}{k}(w^x_0,w^x_{s+\p-k})=0.\tag b$$ 
(Note that if $k>s$ then $(w^x_0,w^x_{s+\p-k})=0$.) This can be         
viewed as an inductive formula for $(w^x_0,w^x_{s+\p})$ (for $k\in[1,2\p+1],k\le s$ we have
$s+\p-k\in[\p,s+\p-1]$). The same inductive formula holds for $|^x_0:^x_{s+\p}|$, see 2.5. It follows that
$(w^x_0,w^x_{s+\p})=|^x_0:^x_{s+\p}|$ for any $s\in\ZZ_{\ge2}$. Hence $(w^x_i,w^x_j)=|^x_i:^x_j|$ for any 
$i,j\in\ZZ$ such that $j-i\ge\p+2$. The last equality also holds for $i,j$ such that $j-i\in[0,\p+1]$ and then by
symmetry, for any $i,j$.

\subhead 2.22\endsubhead
Assume that $\p\ge1$ and that $p_1>\p=p_y=p_x$, $y<x$. Then $a,b,I_\p$ are defined (see 2.17). Assume further that
$b$ is odd. Then 2.19(e) holds for $x$ and also for $y$ instead of $x$; these two identities have the same right 
hand side, hence they have equal left hand sides (after multiplication by $\mu$):
$$\sum_{k\in[0,2\p]}n_kw_k^y=\sum_{k\in[0,2\p]}n_kw_k^x.\tag a$$
Applying $(,w^x_{\p+s})$ (with $s\in\ZZ_{>0}$) to (a) gives
$$(w_0^y,w^x_{\p+s})+\sum_{k\in[1,2\p];k<s}n_k(w_k^y,w^x_{\p+s})=\sum_{k\in[0,2\p];k\le s}n_k(w_k^x,w^x_{\p+s}).
\tag b$$
(Note that if $k\ge s$ then $(w_k^y,w^x_{\p+s})=0$; indeed, we have $0\le k-s\le2\p-1$.)
Here the right hand side is equal to $\sum_{k\in[0,2\p]}n_k|_k^x:^x_{\p+s}|$ by 2.20.
In the left hand side we have for $k\in[1,2\p],k<s$: $(w_k^y,w^x_{\p+s})=(w_0^y,w^x_{\p+s-k})$ and $0<s-k<s$.
Thus (b) can be viewed as an inductive formula for $(w_0^y,w^x_{\p+s})$. This is the same as the inductive
formula 2.6 for $|_0^y:^x_{\p+s}|$. It follows that $(w_0^y,w^x_{\p+s})=|_0^y:^x_{\p+s}|$ for $s\ge1$.
We see that $(w_i^y,w^x_j)=|_i^y:^x_j|$ for any $i,j\in\ZZ$ such that $j-i>\p$.

Applying $(,w^x_{\p-s})$ (with $s\in\ZZ_{\ge0}$) to (a) gives
$$(w^y_{2\p},w^x_{\p-s})+\sum_{k\in[0,2\p-1];k\ge2\p-s}n_k(w_k^y,w^x_{\p-s})
=\sum_{k\in[0,2\p]}n_k(w_k^x,w^x_{\p-s}).\tag c$$
(Note that if $k<2\p-s$ then $(w_k^y,w^x_{\p-s})=0$; indeed we have $0\le k+s\le2\p-1$.)

Here the right hand side is equal to $\sum_{k\in[0,2\p]}n_k|_k^x:^x_{\p-s}|$ by 2.20. In the left hand side we 
have for $k\in[0,2\p-1],k\ge2\p-s$: $(w_k^y,w^x_{\p-s})=(w^y_{2\p},w^x_{\p-s+2\p-k})$ and $0\le s+k-2\p<s$.
Thus (c) can be viewed as an inductive formula for $(w_{2\p}^y,w^x_{\p-s})$. This is the same as the inductive
formula 2.6 for $|_{2\p}^y:^x_{\p-s}|$. It follows that $(w_{2\p}^y,w^x_{\p-s})=|_{2\p}^y:^x_{\p-s}|$ for $s\ge0$.
We see that $(w_i^y,w^x_j)=|_i^y:^x_j|$ for any $i,j\in\ZZ$ such that $i-j\ge\p$.
Since $(w_i^y,w^x_j)=|_i^y:^x_j|=0$ for any $i,j\in\ZZ$ such that $-\p\le i-j<\p$ it follows that
$(w_i^y,w^x_j)=|_i^y:^x_j|=0$ for any $i,j\in\ZZ$.

\subhead 2.23\endsubhead
Assume that $\p\ge1$. Assume further that either $p_1>\p$ and $b$ (see 2.17) is even 
or that $p_1=\p$. We apply 2.16(a),(b) with $c=b+1$ (if $p_1>\p$) or with $c=1$ (if $p_1=\p$). We deduce that
for any $z\in[1,\s]$ such that $p_z=\p$ we have $N^{2\p}w^z_0=\z_zN^{2\p}w^c_0$ where $\z_z\in\kk$. We show that 
$\z_x=\pm1$. If $z=c$ this is clear. We now assume that $z>c$. We have
$$\sum_{k\in[0,2\p]}n_kw^z_k=\z_z\sum_{k\in[0,2\p]}n_kw^c_k.\tag a$$
Taking $(,w^z_\p)$ with (a) we obtain $\z_z(w^c_{2\p},w^z_\p)=2$. In particular $\z_z\ne0$ and 
$2\z_z\i=(w^z_0,w^c_\p)$. Taking $(,w^c_\p)$ with (a) we obtain $(w^z_0,w^c_\p)=2\z_z$. We see that 
$2\z_z\i=2\z_z$ hence $\z_z=\pm1$ as claimed. Replacing $w^z_i$ by $\z_zw^z_i$ for all $z$ such that $p_z=\p$
and all $i\in\ZZ$ we see that we can assume that $\z_z=1$ for all $z$ as above. 

\subhead 2.24\endsubhead
Assume that $\p\ge1$ and that $\p=p_y=p_x$, $y<x$. Assume further that either $p_1>\p$ and $b$ (see 2.17) is even
or that $p_1=\p$. Applying 2.23(a) to $x$ and $y$ instead of $z$ and using that $\z_x=\z_y=1$ we see that
$$\sum_{k\in[0,2\p]}n_kw^x_k=\sum_{k\in[0,2\p]}n_kw^y_k.$$
From this we deduce exactly as in 2.22 that $(w_i^y,w^x_j)=|_i^y:^x_j|=0$ for any $i,j\in\ZZ$.

\subhead 2.25\endsubhead
Assume that $\k=1$ and $\p=1/2$. Assume further that $\s$ is even. We have $x=\s+1$ and $\ps(\s)=-1$. We apply 
2.15(a) with $r=\s$; in this case $\cw$ is spanned by $\{w^t_i;t\in[1,\s],i\in[0,2p_t-1]\}$ and $\cw'$ is spanned 
by $w^x_0$. Using 2.15(a) we see that $g\cw=\cw$ and $g\cw'=\cw'$. Since $g$ is unipotent and $\dim\cw'=1$ we see
that $g=1$ on $\cw'$. Hence $w^x_j=w^x_0$ for all $j\in\ZZ$. Since $\cw$ is $g$-stable, for any 
$r\in[1,\s],i\in\ZZ$ we have $w^r_i\in\cw$; since $(\cw,w^x_0)=0$ and $w^x_j=w^x_0$ we see that for any 
$i,j\in\ZZ$ we have $(w^r_i,w^x_j)=0$. Hence $(w^r_i,w^x_j)=|^r_i:^x_j|$. (See 2.7.) Note also that for any 
$i,j\in\ZZ$ we have $(w^x_i,w^x_j)=(w^x_0,w^x_0)=2$; hence $(w^x_i,w^x_j)=|^x_i:^x_j|$. (See 2.7.)

\subhead 2.26\endsubhead
Assume that $\k=1$ and $\p=1/2$. 
Then $a,b,I_\p$ are defined (see 2.17). We have $b=\s$, $x=\s+1$. Assume further that $\s$ is odd. Then any Jordan
block of $N:V@>>>V$ has size $\ge2$. 
Let $V_1$ be the subspace of $V$ spanned by $\{w^r_i;r\in[a,\s],i\in[0,2p_r-1]\}$.
Let $V_2$ be the subspace of $V$ spanned by $\{w^r_i;r\in[a,\s+1],i\in[0,2p_r-1]\}$.
We have $V_2=V_1\op\kk w^x_0$ with $(V_1,w^x_0)=0$, $(w^x_0,w^x_0)=2$. This, together with the fact that
$(,)_{V_2}$ is nonsingular (see 2.14(b)) implies that $(,)|_{V_1}$ is nonsingular. In particular
the symmetric matrix $((w^r_i,w^{r'}_{i'}))_{r,r'\in[a,\s],i\in[0,2p_r-1],i'\in[0,2p_{r'}-1]}$ is nonsingular.
By the induction hypothesis this is the same as the symmetric matrix
$$(|^r_i:^{r'}_{i'}|)_{r,r'\in[a,\s],i\in[0,2p_r-1],i'\in[0,2p_{r'}-1]}$$ 
which is therefore nonsingular. Since $V_2$ is $g$-stable (see 2.15(a)) we can write 
$$w^a_{2p_a}=\sum_{r\in[a,\s+1],i\in[0,2p_r-1]}C^r_iw^r_i\tag a$$ 
with $C^r_i\in\kk$. For any $r'\in[a,\s],i'\in[0,2p_{r'}-1]$ we have $(w^x_0,w^{r'}_{i'})=0$ hence
$$(w^a_{2p_a},w^{r'}_{i'})=\sum_{r\in[a,\s],i\in[0,2p_r-1]}(w^r_i,w^{r'}_{i'})C^r_i$$
that is,
$$|^a_{2p_a}:^{r'}_{i'}|=\sum_{r\in[a,\s],i\in[0,2p_r-1]}|^r_i:w^{r'}_{i'}|C^r_i.$$
This can be regarded as a system of linear equations with unknowns $C^r_i$ and with a nonsingular matrix. This 
is the same as the system of linear equations defining $c^r_i$ in 2.8. It follows that $C^r_i=c^r_i$ for 
$r\in[a,\s],i\in[0,2p_r-1]$.

Assume that $C^x_0=0$. Then from (a) we see that $w^a_{2p_a}\in V_1$. 
Note that $V_1$ is spanned by $\{N^iw^r_0;r\in[a,\s],i\in[0,2p_r-1]\}$.
Since $w^a_{2p_a}$ is equal to $N^{2p_a}w^a_0$ plus a linear combination of elements in $V_1$ it follows that
$N^{2p_a}w^a_0\in V_1$. We show that $NV_1\sub V_1$. It is enough to show that $NN^iw^r_0\in V_1$ for any
$r\in[a,\s],i\in[0,2p_r-1]$. If $i\in[0,2p_r-2]$ this is obvious. If $i=2p_r-1$ then 
$NN^iw^r_0=N^{2p_r}w^r_0$ which by 2.16(a) belongs to the space generated by $V_1$ and by 
$N^{2p_a}w^a_0$ (which also belongs to $V_1$) hence $NN^iw^r_0\in V_1$. Thus $NV_1\sub V_1$ and $gV_1\sub V_1$.
Using the decomposition $V_2=V_1\op\kk w^x_0$ in which $\kk w^x_0$ is the perpendicular to $V_1$ in $V_2$ it
follows that the line $\kk w^x_0$ is also $g$-stable. Since $g$ is unipotent it must act on $\kk w^x_0$ as
identity. We see that $N:V_2@>>>V_2$ has at least one Jordan block of size $1$.
Using 2.15(a) it follows that $N:V@>>>V$ has at least one Jordan block of size $1$. But this is not the case.
We have therefore proved that $C^x_0\ne0$. Taking self inner products in (a) and using the induction hypothesis
we obtain
$$0=\sum_{r,r'\in[a,\s];i\in[0,2p_r-1];i'\in[0,2p_{r'}-1]}c^r_ic^{r'}_{i'}|^r_i:^{r'}_{i'}|+2(C^x_0)^2.$$
(Note that $|^a_{2p_a}:^a_{2p_a}|=0$.) Hence
$$C^x_0=\e_x\sqrt{\nu/2}$$
where
$$\nu=-\sum_{r,r'\in[a,\s];i\in[0,2p_r-1];i'\in[0,2p_{r'}-1]}c^r_ic^{r'}_{i'}|^r_i:^{r'}_{i'}|$$
and $\e_x=\pm1$. Replacing $w^x_i$ by $\e_xw^x_i$ for any $i\in\ZZ$ we see that we can assume that $\e_x=1$.
Comparing with 2.8 we see that $C^x_0=c^x_0$. In particular, $c^x_0\ne0$. We have 
$$w^x_0=(c^x_0)\i(w^a_{2p_a}-\sum_{r\in[a,\s],i\in[0,2p_r-1]}c^r_iw^r_i).$$
Applying $g^h$ ($h\in\ZZ)$ we obtain
$$w^x_h=(c^x_0)\i(w^a_{2p_a+h}-\sum_{r\in[a,\s],i\in[0,2p_r-1]}c^r_iw^r_{i+h}).\tag b$$
For any $r'\in[a,\s],i'\in[0,2p_{r'}-1]$ we have (using (b) and the induction hypothesis):
$$(w^{r'}_{i'},w^x_h)=
(c^x_0)\i(|^a_{2p_a+h}:^{r'}_{i'}|-\sum_{r\in[a,\s],i\in[0,2p_r-1]}c^r_i|^r_{i+h}:^{r'}_{i'}|)$$
hence 
$$(w^x_h,w^{r'}_{i'})=|^x_h:^{r'}_{i'}|.$$
(See 2.8.) For $h,h'\in\ZZ$ we have (using (b) and the induction hypothesis):
$$\align&(w^x_h,w^x_{h'})=(c^x_0)^{-2}(|^a_{2p_a+h}:^a_{2p_a+h'}|\\&
-\sum_{r\in[a,\s],i\in[0,2p_r-1]}c^r_i(|^r_{i+h}:^a_{2p_a+h'}|+|^r_{i+h'}:^a_{2p_a+h}|\\&
+\sum_{r,r'\in[a,\s],i\in[0,2p_r-1],i'\in[0,2p_{r'}-1]}c^r_ic^{r'}_{i'}|^r_{i+h}:^{r'}_{i'+h'}|).\endalign$$
Hence
$$(w^x_h,w^x_{h'})=|^x_h:^x_{h'}|.$$
(See 2.8.) This completes the inductive proof of Proposition 2.13.

\head 3. Proof of Theorem 0.2\endhead
\subhead 3.1\endsubhead
Let $V,Q,(,),\nn,n,\k,Is(V)$ be as in 1.1. Let $Is(V)^0$ be the identity component of $Is(V)$. A subspace $V'$ of
$V$ is said to be isotropic if $(,)$ and $Q$ are zero on $V'$. Let $\cf$ be the set of all sequences 
$V_*=(0=V_0\sub V_1\sub V_2\sub\do\sub V_\nn=V)$ of subspaces of $V$ such that $\dim V_i=i$ for $i\in[0,\nn]$, 
$Q|_{V_i}=0$ and $V_i^\pe=V_{\nn-i}$ for all $i\in[0,n]$. (For such $V_*$, $V_i$ is an isotropic subspace for 
$i\in[0,n]$). Now $Is(V)$ acts naturally (transitively) on $\cf$. 

Let $p_1\ge p_2\ge\do\ge p_\s$ (or $p_*$) be as in 1.2. (If $\k=0,Q\ne0$ we assume that $\s$ is even.)
Let $(V_*,V'_*)\in\cf\T\cf$ be such that for any $r\in[1,\s]$ we have
$$\dim(V'_{p_{<r}+i}\cap V_{p_{<r}+i})=p_{<r}+i-r,\qua \dim(V'_{p_{<r}+i}\cap V_{p_{<r}+i+1})=p_{<r}+i-r+1$$
if $i\in[1,p_r-1]$;
$$\dim(V'_{p_{\le r}}\cap V_{\nn-p_{<r}-1})=p_{\le r}-r,\qua \dim(V'_{p_{\le r}}\cap V_{\nn-p_{<r}})
=p_{\le r}-r+1.$$
Here $p_{\le r}=\sum_{i\in[1,r]}p_i,p_{<r}=\sum_{i\in[1,r-1]}p_i$.
Let $g\in Is(V)^0,\tg\in Is(V)^0$ be unipotent elements such that $gV_*=V'_*$, $\tg V_*=V'_*$. Let $N=g-1$, 
$\tN=\tg-1$. We assume that 

if $Q=0$ or $p=2$ then $\cm(N,V)=\cm(\tN,V)$ consists of $2p_1\ge2p_2\ge\do\ge2p_\s$ (and $1$ if $\k=1$);
 
if $Q\ne0$ and $p\ne2$ then $\cm(N,V)=\cm(\tN,V)$ consists of $2p_1+\ps(1)\ge2p_2+\ps(2)\ge\do\ge2p_\s+\ps(s)$ 
(and $1$ if $\k=1$)
\nl
with $\ps(r)$ as in 2.1.

We show:

(a) {\it there exists $T\in Is(V)^0$ such that $\tg=TgT\i$, $T(V_*)=V_*$, $T(V'_*)=V'_*$.}
\nl
Let $v_1,v_2,\do,v_{\s+\k}$ be the sequence of vectors associated in \cite{\WE, 3.3} to $V_*,V'_*,g$ (each $v_i$
is uniquely defined up to multiplication by $\pm1$). Let $\tv_1,\tv_2,\do,\tv_{\s+\k}$ be the analogous sequence 
defined in terms of $V_*,V'_*,\tg$ instead of $V_*,V'_*,g$. From the definitions, $w^r_i:=g^{-p_r+i}v_r$ 
$(r\in[1,\s+\k],i\in\ZZ)$ is a $(g,p_*)$-adapted collection of vectors in $V$. Similarly,
$\tw^r_i:=\tg^{-p_r+i}\tv_r$ $(r\in[1,\s+\k],i\in\ZZ)$ is a $(\tg,p_*)$-adapted collection of vectors in $V$. 
Using 1.3(a) we see that there is a unique linear isomorphism $T:V@>>>V$ such that $T(w^r_i)=\tw^r_i$ for all 
$r\in[1,\s+\k],i\in[0,2p_r-1]$.

Using Proposition 1.6 (in the case where $Q=0$ or $p=2$) or Proposition 2.13 (in the case where $Q\ne0$ and 
$p\ne2$) we see that there exist functions $\e,\ti\e:[1,\s+\k]@>>>\{1,-1\}$, $t\m\e_t,t\m\ti\e_t$, such that 
$(\e_yw^y_i,\e_xw^x_j)=(\ti\e_y\tw^y_i,\ti\e_x\tw^x_j)$ for all $x,y\in[1,\s+\k],i,j\in\ZZ$. Replacing $v_r$ by 
$\e_rv_r$ and $\tv_r$ by $\ti\e_r\tv_r$ we see that we can assume that
$$(w^y_i,w^x_j)=(\tw^y_i,\tw^x_j)\tag b$$
for all $x,y\in[1,\s+\k],i,j\in\ZZ$. Thus we have $(T(w^y_i),T(w^x_j))=(w^y_i,w^x_j)$ for all 
$x,y\in[1,\s+\k],i\in[0,2p_x-1],j\in[0,2p_y-1]$. Moreover, if $Q\ne0,p=2$, we have 
$Q(w^y_i)=Q(\tw^y_i)=Q(T(w^y_i))=0$ for $y\in[1,r],i\in[0,2p_y-1]$ and 
$Q(w^{\s+1}_0)=Q(\tw^{\s+1}_0)=Q(T(w^{\s+1}_0))=1$ if $\k=1$. Since 
$\{w^y_i;y\in[1,\s+\k],i\in[0,2p_y-1]\}$ is a basis of $V$, we see that $T\in Is(V)$. 

Let $x,y\in[1,\s+\k],i\in\ZZ,j\in[0,2p_y-1]$. Since $T\in Is(V)$ we have $(T(w^y_i),T(w^x_j))=(w^y_i,w^x_j)$ that
is $(T(w^y_i),\tw^x_j)=(w^y_i,w^x_j)$. Moreover we have $(\tw^y_i,\tw^x_j)=(w^y_i,w^x_j)$, see (b). Thus 
$(T(w^y_i),\tw^x_j)=(\tw^y_i,\tw^x_j)$ and $(T(w^y_i)-\tw^y_i,\tw^x_j)=0$. Since the elements $\tw^x_j$, 
$(x\in[1,\s+\k],j\in[0,2p_x-1])$ form a basis of $V$ it follows that $\x:=T(w^y_i)-\tw^y_i\in V^\pe$. If $\k=0$ or
$p\ne2$ it folows that $\x=0$. If $\k=1$ and $p=2$ we have $T(w^y_i)=\tw^y_i+\x$ hence 
$Q(T(w^y_i))=Q(\tw^y_i)+Q(\x)$. Since $T\in Is(V)$ this implies $Q(w^y_i)=Q(\tw^y_i)+Q(\x)$. But from the 
definitions we have $Q(w^y_i)=Q(\tw^y_i)$ hence $Q(\x)=0$ so that $\x=0$. We see that in any case $\x=0$ that is
$T(w^y_i)=\tw^y_i$ (for $y\in[1,\s+\k],i\in\ZZ$).

For $y\in[1,\s+\k],i\in\ZZ$ we have $T(g(w^y_i))=T(w^y_{i+1})=\tw^y_{i+1}=\tg(\tw^y_i)=\tg(T(w^y_i))$. Since the 
elements $w^y_i$ $(y\in[1,\s+\k],i\in[0,2p_y-1])$ form a basis of $V$, it follows that $Tg=\tg T$.

From the definitions, for any $r\in[1,\s],i\in[0,p_r]$, the subspace $V_{p_{<r}+i}$ is generated by $w^t_h$ 
$(t<r,h\in[p_t,2p_t-1])$ and by $w^r_h$ $(h\in[p_t,p_t+i-1])$; similarly, the subspace $V_{p_{<r}+i}$ is 
generated by $\tw^t_h$ $(t<r,h\in[p_t,2p_t-1])$ and by $\tw^r_h$ $(h\in[p_t,p_t+i-1])$. Applying $T$ we see that 
the subspace $T(V_{p_{<r}+i})$ is generated by $T(w^t_h)=\tw^t_h$ $(t<r,h\in[p_t,2p_t-1])$ and by 
$T(w^r_h)=\tw^r_h$ $(h\in[p_t,p_t+i-1])$. It follows that $T(V_{p_{<r}+i})=V_{p_{<r}+i}$. Thus we have 
$T(V_*)=V_*$. From $gV_*=V'_*$, $\tg V_*=V'_*$, $T(V_*)=V_*$, $Tg=\tg T$ we deduce $T(V'_*)=V'_*$.

We show that $T\in Is(V)^0$. If $Q=0$ or $k=1$ this is obvious. Now assume that $Q\ne0,\k=0$. We have 
$T(V_n)=V'_n=g(V_n)$. Since $g\in Is(V)^0$, we see that $V_n$ and $g(V_n)$ are in the same $Is(V)^0$-orbit. Hence
$T(V_n)$ and $V_n$ are in the same $Is(V)^0$-orbit that is, $T\in Is(V)^0$. This completes the proof of (a).

\subhead 3.2\endsubhead
We now prove Theorem 0.2. 
Using \cite{\WE, 5.2(a)} we see that it is enough to prove the theorem for a particular $w\in C_{min}$; then it
will automatically hold for any $w\in C_{min}$. We can assume that $G$ is adjoint and almost simple. Moreover for 
each isogeny class of almost simple groups it is enough to consider one group in the isogeny class and the result
will be automatically true for the other groups in the isogeny class. If $G$ is of type $A$ the result is easily
proved; we omit the details. 
If $G$ is of type $B,C$ or $D$, we can assume that $G=Is(V)^0$ with $V$ as in 1.1. Then for any $w\in C_{min}$
we have $\fB^\g_w\ne\em$ by \cite{\WE, 4.6(a)} and for a specific $w\in C_{min}$ the $G$ action on $\fB^\g_w$ is 
transitive by 3.1(a). 
 Now assume that $G$ is simple of exceptional type. We can assume that $\kk$ is an algebraic closure 
of a finite field $\FF_q$ with $q$ elements. (We assume also that $q-1$ is sufficiently
divisible.) We choose an $\FF_q$-split rational structure on $G$ with Frobenius map
$F:G@>>>G$. Note that $F(\g)=\g$. Now $F$ induces a morphism $\cb@>>>\cb$ denoted again by $F$ and $\fB^\g_w$ has
a natural Frobenius map $(g,B)\m(F(g),F(B))$ denoted again by $F$. We calculate the number of fixed points of 
$F:\fB^\g_w@>>>\fB^\g_w$: using the method of \cite{\WE, 1.2} this is reduced to a computer calculation which 
shows that this number is equal to $|G(\FF_q)|$. Since this holds when $q$ is replaced by any power of $q$
we conclude that $\fB^\g_w$ is an irreducible variety of dimension equal to $\dim(G)$. By \cite{\WE, 5.2} the
$G$-action on $\fB^\g_w$ has finite isotropy groups. Hence each $G$-orbit on $\fB^\g_w$ has dimension equal to 
$\dim(G)$. It follows that $\fB^\g_w$ is a single $G$-orbit. This completes the proof of Theorem 0.2.

\subhead 3.3\endsubhead
Here is a complement to Theorem 0.2. Assume that either $p$ is not a bad prime for $G$ or that $G$ is a simple 
exceptional group with $p$ arbitrary. Let $C\in\uWW_{el}$ and let $w\in C_{min}$. 
Let $\g$ be a unipotent class in class in $G$ such that $\g\ne\Ph(C)$. Then

(a) {\it $\fB^\g_w$ is a union of infinitely many $G$-orbits.}
\nl
If $p$ is not a bad prime for $G$ then this follows from \cite{\WE, 5.8(b)}. Now assume that $G$ is a simple 
exceptional group. As in 3.2 we calculate the number of fixed points of $F:\fB^\g_w@>>>\fB^\g_w$ with the aid of 
a computer ($F$ is relative to an $\FF_q$-structure as in 3.2.) We find that this number is either $0$ (for all 
$q$) or a polynomial in $q$ of degree $>\dim(G)$. The conclusion follows. We expect that (a) holds without 
assumptions on $G$ or $p$.

\head 4. Not necessarily elliptic elements\endhead
\subhead 4.1\endsubhead
Let $C\in\uWW$. Let $\g=\Ph(C)$, see 0.1. Let $w\in C_{min}$. We consider the following property:

(a) {\it $\fB^\g_w$ is a single $G$-orbit and $\dim\fB^\g_w=(1/2)(\dim G_{ad}+\dim\g-\mu(w)+\ul(w))$.}
\nl
We show that (a) holds if $C\in\uWW_{el}$. The first assertion of (a) follows from 0.2. In particular we have
$\fB^\g_w\ne\em$ so that $C\dsv\g$. Using \cite{\WE, 4.4(b)} and \cite{\WE, 5.7(iv)} we see that 
$\dim\fB^\g_w=\dim G_{ad}$ and $\dim G_{ad}-\dim\g=\ul(w)$. Note that $\mu(w)=0$. Hence 
$(1/2)(\dim G_{ad}+\dim\g-\mu(w)+\ul(w))=\dim G_{ad}$ and (a) follows in this case.

We now drop the assumption that $C\in\uWW_{el}$.

\proclaim{Proposition 4.2}(i) If 4.1(a) holds for some $w\in C_{min}$ then it holds for any $w\in C_{min}$.

(ii) Assume that $G$ is simple of type $\ne B,C,D$; then 4.1(a) holds for any $w\in C_{min}$.
\endproclaim
The proof is given in 4.5, 4.6. We expect that (ii) holds without any assumption.

\subhead 4.3\endsubhead
Let $C\in\uWW$. Let $w\in C_{min}$. Let $K$ be the set of all elements of $S$ which appear in some/any reduced
expression of $w$. Let $\cc$ be the conjugacy class of $w$ in $\WW_K$. From \cite{\GP, 3.1.12} we see that $w$ 
has minimal length among the elements of $\cc$. Let $X$ be the 
variety consisting of the pairs $(P,L)$ where $P\in\cp_K$ and $L$ is a Levi subgroup of $P$. For any $(P,L)\in X$
we may identify canonically $\WW_K$ with the Weyl group of $L$ and we denote by $\g_{P,L}$ the unipotent 
conjugacy class in $L$ associated to $\cc$ by the analogue of $\Ph$ for $L$. Note that 
$x\g_{P,L}x\i=\g_{xPx\i,xLx\i}$ for any $x\in G$. From \cite{\WE, 1.1} have $\g_{P,L}\sub\g$.
Let $Y$ be the variety consisting of all triples $(g,B,L)$ where $B\in\cb$, $L$ is a Levi subgroup of $P^K_B$,
$g\in\g_{P^K_B,L}$ and $(B,gBg\i)\in\co_w$. Note that $G$ acts on $Y$ by conjugation on all factors. We have the 
following result:

(a) {\it $Y$ is a single $G$-orbit; if $\ci$ is the isotropy group of $(g,B,L)\in Y$ then 
$\ci^0=\cz_L^0$ (a torus of dimension $\dim\cz_G+\mu(w))$.}
\nl
Note that $G$ acts on $X$ transitively (by conjugation) and the isotropy group of $(P,L)$ is $L$. Define 
$r:Y@>>>X$ by $(g,B,L)\m(P^K_B,L)$. This map is $G$-equivariant. Hence it is enough to show that for any 
$(P,L)\in X$, $r\i(P,L)$ is a single $L$-orbit and the isotropy group in $L$ at $(g,B,L)$ has 
identity component equal to $\cz_L^0$. But $r\i(P,L)$ can be identified with $\fB^{\g'}_w$ (defined like
$\fB^\g_w$ in terms of $L,\g'=\g_{P,L},w$ instead of $G,\g,w$). Then the desired statement follows 
from 0.2 applied to $L$ and from \cite{\WE, 5.2} applied to $L$.

Let $Z$ be the variety consisting of all $(g,B)\in\fB^\g_w$ such that for some Levi subgroup $L$ of $P^K_B$ we 
have $g\in\g_{P^K_B,L}$ (hence $g\in L$). Note that $G$ acts on $Z$ by conjugation on both factors. Define 
$\p:Y@>>>Z$ by $\p(g,B,L)=(g,B)$. This map is clearly surjective; it is also $G$-equivariant. Hence using (a) we
see that $Z$ is a single $G$-orbit. The fibre of $\p$ at $(g,B)\in Z$ can be identified with the set
$R$ of Levi subgroups $L'$ of $P:=P^K_B$ such that $g\in\g_{P,L'}$ (hence $g\in L'$). Let $L\in R$ and let $U$ be
the unipotent radical of $P$. Then 
$$R=\{uLu\i;u\in U,g\in uLu\i\}=\{uLu\i; u\in U,ug=gu\}\cong R_0$$
where $R_0=\{u\in U;gu=ug\}$. By \cite{\HS}, $R_0$ is a connected unipotent group and by \cite{\IC, 2.9(a)}, we 
have $\dim(R_0)=\dim U-(1/2)(\dim\g-\dim\g')$ where $\g'=\g_{P,L}$. We apply \cite{\WE, 5.7(iv)} to $L$ 
(using 0.2 and \cite{\WE, 4.4(b)} for $L$); we obtain 
$\dim\g'=\dim L-\dim\cz_L-\ul(w)$. Note also that $\dim U=(\dim G-\dim L)/2$ and $\dim\cz_L=\dim\cz_G+\mu(w)$. 
Hence 
$$\dim(R_0)=(1/2)(\dim G_{ad}-\dim\g-\mu(w)-\ul(w)).$$

\subhead 4.4\endsubhead
In the setup of 4.3 we assume that $\kk$ is an algebraic closure of a finite field $\FF_q$ with $q$ elements. (We 
assume also that $q-1$ is sufficiently divisible.) We choose an $\FF_q$-split rational structure on $G$ with 
Frobenius map $F:G@>>>G$. Note that $F(\g)=\g$. Now $F$ induces morphisms $X@>>>X,Y@>>>Y,Z@>>>Z$,
$\fB^\g_w@>>>\fB^\g_w$ denoted again by $F$. We show that
$$|Z^F|=|G_{ad}^F|q^{-(1/2)(\dim G_{ad}-\dim\g-\mu(w)-\ul(w))}(q-1)^{-\mu(w)}.\tag a$$
The map $\p:Y@>>>Z$ restricts to a map $\p':Y^F@>>>Z^F$. Let $(g,B)\in Z^F$. By the arguments in 
4.3, $\p'{}\i(g,B)$ is the set of fixed points of $F$ on a principal homogeneous space of a connected
(unipotent) group hence this set is nonempty and has cardinal equal to $q$ raided to a power equal to the 
dimension of that unipotent group. We see that $|Z^F|=|Y^F|q^{-(1/2)(\dim G_{ad}-\dim\g-\mu(w)-\ul(w))}$. It 
remains to show that $|Y^F|=|G_{ad}^F|(q-1)^{-\mu(w)}$. The map $r:Y@>>>X$ restricts to a map set 
$r':Y^F@>>>X^F$. Since $|X^F|=|G^F|/|L^F|$ and $|L^F|=|L_{ad}^F|(q-1)^{\mu(w)}$) it is enough to show that for any
$(P,L)\in X^F$ we have $|r'{}\i(P,L)|=|L_{ad}^F|$. This follows from the fact that $r'{}\i(P,L)=(r\i(P,L))^F$ and
$r\i(P,L)$ is a homogeneous space for $L_{ad}$ with finite isotropy groups (see 4.1). This completes the proof of
(a). 

\subhead 4.5\endsubhead
We prove 4.2(i). We can assume that we are in the setup of 4.4.
Assume that 4.1(a) holds for some $w\in C_{min}$. Since $Z\sub\fB^\g_w$ and $Z$ is a single $G$-orbit it
follows that $\fB^\g_w=Z$ and $|(\fB^\g_w)^F|$ is given by the right hand side of 4.4(a). Since
$|(\fB^\g_w)^F|$ is independent of the choice of $w$ in $C_{min}$ we deduce that if $w'\in C_{min}$ then
$|(\fB^\g_{w'})^F|$ is given by the right hand side of 4.4(a). Let us define $Z'$ in terms of $w'$ in the same way
as $Z$ was defined in terms of $w$. Then 4.4(a) applied to $w'$ instead of $w$ shows that $|Z'{}^F|$ is again 
given by the right hand side of (a). (Note that $\ul(w)=\ul(w'),\mu(w)=\mu(w')$.) Thus
$|(\fB^\g_{w'})^F|=|Z'{}^F|$. Since $Z'{}^F\sub(\fB^\g_{w'})^F$ it follows that $(\fB^\g_{w'}-Z')^F=\em$. Since 
this holds when $F$ is replaced by any power of $F$ it follows that $\fB^\g_{w'}-Z'=\em$ that is $\fB^\g_{w'}=Z'$.
Since $Z'$ is a single $G$-orbit it follows that $\fB^\g_{w'}$ is a single $G$-orbit. Since
$\dim Z'=(1/2)(\dim G_{ad}+\dim\g-\mu(w')+\ul(w'))$ (by the arguments in 4.3) we see that
$\dim\fB^\g_{w'}=(1/2)(\dim G_{ad}+\dim\g-\mu(w')+\ul(w'))$. This proves 4.2(i).

\subhead 4.6\endsubhead
We prove 4.2(ii). It is easy to see that 4.2(ii) holds if $G$ is of type $A$. We now assume that $G$ is of 
exceptional type. We can assume that we are in the setup of 4.4.
We compute the number $|(\fB^\g_w)^F|$; using the method of \cite{\WE, 1.2} this is reduced to a 
computer calculation which shows that this is number is equal to the right hand side of 4.4(a).
It follows that $|(\fB^\g_w)^F|=|Z^F|$. Since $Z^F\sub(\fB^\g_w)^F|$ it follows that $(\fB^\g_w-Z)^F=\em$. Since 
this holds when $F$ is replaced by any power of $F$ it follows that $\fB^\g_w-Z=\em$ that is $\fB^\g_w=Z$.
Since $Z$ is a single $G$-orbit it follows that $\fB^\g_w$ is a single $G$-orbit. Since
$\dim Z=(1/2)(\dim G_{ad}+\dim\g-\mu(w)+\ul(w))$ (by the arguments in 4.3) we see that
$\dim\fB^\g_w=(1/2)(\dim G_{ad}+\dim\g-\mu(w)+\ul(w))$. This proves 4.2(ii).

\subhead 4.7\endsubhead
Let $C\in\uWW$ and let $w\in C_{min}$. From the proof in 4.5 we see that property 4.1(a) holds if and only if
the following holds:

(a) {\it for any $(g,B)\in\fB^\g_w$ there exists a Levi subgroup $L$ of $P^K_B$ such that $g\in\g_{P^K_B,L}$ 
(hence $g\in L$). (Notation of 4.3).}

\subhead 4.8\endsubhead
We give an alternative definition of the map $\Ph:\uWW@>>>\uuG$ of \cite{\WE, 4.5} assuming that $G$ is simple of
exceptional type.

(a) {\it Let $C\in\uWW$. There exists a unique $\g\in\uuG$ such that $C\dsv\g$ and such that if $\g'\in\uuG$ and 
$C\dsv\g'$ then $\g\sub\bar\g'$. We have $\g=\Ph(C)$.}
\nl
Note that when $p$ is not a bad prime this is already in \cite{\WE}. The proof in the case where $p$ is a
bad prime is similar (based on computer calculation). We expect that (a) holds for arbitrary $G$.

\enddocument